\def\docdate{May 12, 2006}

\documentclass[12pt,reqno]{amsart}

\title[Closure properties and heavy tails]{Closure properties and heavy tails: \\random vectors in presence of dependence}

\author{Dimitrios G. Konstantinides,\quad Charalampos  D. Passalidis} 

\date{\docdate}

\usepackage{amssymb}
\usepackage{amsmath}
\usepackage{amsthm}
\usepackage{latexsym}

\newtheorem{theorem}{Theorem}[section]
\newtheorem{lemma}{Lemma}[section]
\newtheorem{corollary}{Corollary}[section]
\newtheorem{remark}{Remark}[section]
\newtheorem{definition}{Definition}[section]
\newtheorem{proposition}{Proposition}[section]
\newtheorem{example}{Example}[section]
\newtheorem{assumption}{Assumption}[section]
\numberwithin{equation}{section}
\newcommand{\bth}{\begin{theorem}}
\newcommand{\ethe}{\end{theorem}}

\newcommand{\bre}{\begin{remark}}
\newcommand{\ere}{\end{remark}}

\newcommand{\ble}{\begin{lemma}}
\newcommand{\ele}{\end{lemma}}
\newcommand{\bde}{\begin{definition}}
\newcommand{\ede}{\end{definition}}

\newcommand{\bco}{\begin{corollary}}
\newcommand{\eco}{\end{corollary}}

\newcommand{\bpr}{\begin{proposition}}
\newcommand{\epr}{\end{proposition}}

\newcommand{\bexer}{\begin{exercise}}
\newcommand{\eexer}{\end{exercise}}
\newcommand{\breh}{\begin{hint}}
\newcommand{\ereh}{\end{hint}}

\newcommand{\halmos}{\hfill \qed}

\newcommand{\bexam}{\begin{example}}
\newcommand{\eexam}{\end{example}}

\newcommand{\pr} {{\bf Proof.}}

\newcommand{\bfi}{\begin{fig}}
\newcommand{\efi}{\end{fig}}

\newcommand{\beao}{\begin{eqnarray*}}
\newcommand{\eeao}{\end{eqnarray*}\noindent}

\newcommand{\beam}{\begin{eqnarray}}
\newcommand{\eeam}{\end{eqnarray}\noindent}
\newcommand{\E}{\mathbf{E}}
\newcommand{\PP}{\mathbf{P}}

\newcommand{\xto}{x\to\infty}

\newcommand{\bH}{\overline{H}}
\newcommand{\bF}{\overline{F}}

\newcommand{\bG}{\overline{G}}

\newcommand{\bbb}{{\mathbb B}}
\newcommand{\bbz}{{\mathbb Z}}
\newcommand{\bbn}{{\mathbb N}}

\newcommand{\vep}{\varepsilon}

\begin{document}
\date{{\small \today}}

\begin{abstract}
This paper is related to closure properties of heavy-tailed random variables and vectors. Firstly, we consider two random variables $X$ and $Y$ with distributions $F$ and $G$ respectively. We assume that the random variables satisfy a type of a weak dependence structure. Under mild conditions, we study if the distribution of $X\,Y$, denoted by $H$, belongs to the same distribution class of $F$. The closure property  under this specific weak dependence structure, in the classes $\mathcal{C},\mathcal{D}, \mathcal{M}, \mathcal{M^*}, \mathcal{OS}, \mathcal{OL}, \mathcal{P_D}$ and $\mathcal{K}$ is established. Secondly, we introduce a distribution class, which satisfies some closure properties with respect to product and finite mixture. Further, we provide some applications on randomly weighted sums and on discrete-time risk model with dependent actuarial and financial risks. The multivariate regular variation is already well-established, but in the case of other heavy tailed classes, the notions of their multivariate distributions is missing. For this reason we introduce the class of dominatedly varying and positively decreasing random vectors and we study the closure property of the independent scalar product. Furthermore, we explore the closure property of the first class under sum and mixture, and we examine the distribution of randomly stopped sums, whose summands are random vectors that belongs to this class. Some of these properties hold and for positively decreasing random vectors.
\end{abstract}

\maketitle



\section{Introduction} \label{sec.KP.1}

Let $X$ and $Y$ be random variables with distribution functions $F$ and $G$ respectively, with support $[0,\,\infty)$. In some case the support of $F$ can be extended on the whole real axis. The distribution of the product $X\,Y$ is denoted by $H(x):=\PP[X\,Y \leq x]$. 

Such products have been studied recently, because of several applications in actuarial and financial mathematics, especially when the distribution of one factor has heavy tails. For example, in risk theory the random variable $X$ depicts the main loss (or gain) in some portfolio over a certain time interval, while $Y$ corresponds to the discount factor, whence the product $X\,Y$ represents the unknown discounted future loss (or gain). For some applications on risk theory and risk management with product distributions see \cite{tang:2006}, \cite{xu:cheng:wang:chang:2018} and \cite{tang:yuan:2014}. For some applications of such products on sample path of infinitely divisible stochastic processes and stochastic recurrence equations, see in \cite{rosinski:samorodnitsky:1993} and \cite{konstantinides:mikosch:2005}.

Next, we denote by $\bF:=1-F$ the distribution tail, hence $\bF(x)=\PP[X>x]$ and assume $\bF(x)>0$ for any $x \geq 0$. For some random variable $X$ with distribution $F$, we denote by $S_X$ (or $S_F$) the support of $F$. For two positive functions $f(x)$ and $g(x)$, the asymptotic relation $f(x)=o[g(x)]$, as $\xto$, denotes
\beao
\lim_{\xto}\dfrac{f(x)}{g(x)} = 0\,,
\eeao
the asymptotic relation $f(x)=O[g(x)]$, as $\xto$, means that
\beao
\limsup_{\xto} \dfrac{f(x)}{g(x)} < \infty\,,
\eeao
and the asymptotic relation $f(x)\asymp g(x)$, as $\xto$, holds if  both $f(x)=O[g(x)]$ and $g(x)=O[f(x)]$. For two real number $x,\,y$, we denote $x^{+}:=\max\{x,0\}$, $x\vee y :=\max \{x,\,y\}$, $x\wedge y :=\min\{x,\,y\}$. For two $d$-dimensional random vectors $\bf b$ and $\bf t$, their Hadamard product is defined as the vector ${\bf b}\,{\bf t}:=(b_1\,t_1,\,\ldots,\,b_d\,t_d)$. 

Let us introduce some heavy tailed distribution classes. We say that a distribution $F$ is heavy tailed and we write $F \in \mathcal{K}$ if holds
\beao
\int_{-\infty}^{\infty} e^{\varepsilon\,y}\,F(dy) = \infty\,,
\eeao
for all $\varepsilon >0$. Further, for a distribution $F$ whose support is the whole real axis, we say that it has long tail and we denote $F \in \mathcal{L}$, if for all (or equivalently, for some) $a > 0$ it holds
\beao
\lim_{\xto} \dfrac{\bF(x-a)}{\bF(x)}=1\,.
\eeao 
The following distribution class was introduced in \cite{chistyakov:1964}. For a distribution $F$ with support $[0,\,\infty)$, we say that it is subexponential and we denote $F \in \mathcal{S}$, if for all (or equivalently, for some) integer $n\geq 2$ it holds
\beao
\lim_{\xto}\dfrac{\overline{F^{n*}}(x)}{\bF(x)} = n\,.
\eeao
Now for a distribution $F$ with a support the whole real axis $(-\infty,\,\infty)$, we say that is subexponential  and we write $F \in \mathcal{S}$, if $F^+ \in \mathcal{S}$, where $F^+(x) := F(x)\,{\bf 1}_{\{x\geq 0\}}$, with ${\bf 1}_A$ the indicator function of the event $A$ (see in \cite{pakes:2004}). For applications of subexponential distributions in actuarial and financial mathematics see for example in \cite{asmussen:albrecher:2010}, \cite{rolski:schmidli:schmidt:teugels:1999}, \cite{schmidli:2017} and \cite{Foss:Korshunov:Zachary:2013}.

The following distribution class was introduced in \cite{feller:1969}. We say that the distribution $F$ belongs to the class of dominatedly varying distributions, symbolically $F \in \mathcal{D}$, if for all (or equivalently, for some) $b \in (0,\,1)$ it holds
\beao
\limsup_{\xto} \dfrac{\bF(bx)}{\bF(x)} < \infty\,.
\eeao
It is known that $\mathcal{D} \nsubseteq \mathcal{S}$ and $\mathcal{S} \nsubseteq \mathcal{D}$, but $\mathcal{D}\cap \mathcal{S}\equiv \mathcal{D}\cap \mathcal{L} \neq \emptyset$, see \cite{goldie:1978}.

Another smaller class is the class of consistently varying distributions, symbolically $\mathcal{C}$. We say that a distribution $F$ has consistently varying tail, if:
\beao
	\lim_{v\downarrow 1}\liminf_{\xto}\dfrac{\overline{F}(v\,x)}{\overline{F}(x)}=1\,,\quad \text{or, equivalently} \qquad 
	\lim_{b\uparrow 1}\limsup_{\xto}\dfrac{\overline{F}(b\,x)}{\overline{F}(x)}=1\,.
\eeao

We say that a distribution $F$ is regularly varying with index $\alpha >0$, symbolically $F \in \mathcal{R}_{-\alpha}$, if for any $t>0$, it holds
\beao
\lim_{\xto} \dfrac{\bF(t\,x)}{\bF(x)} = t^{-\alpha}\,.
\eeao
The relations $\bigcup_{\alpha>0} \mathcal{R}_{-\alpha} \subsetneq \bigcup_{0\leq\alpha\leq\beta}\mathcal{ERV}(-\alpha,-\beta) \subsetneq\mathcal{C}\subsetneq \mathcal{D} \cap \mathcal{L} \subsetneq \mathcal{L} \subsetneq \mathcal{K}$ can be found in  \cite{Bingham:Goldie:Teugels:1987}, \cite{Konstantinides:2018}.

A distribution class, that contains also light tailed distributions, is the class of positively decreasing distributions, suggested in \cite[Unpublished appendix]{dehaan:resnick:1984}. We say that a distribution $F$ is  positively decreasing, symbolically $F \in \mathcal{P_D}$, if for all (or equivalently, for some) $v>1$ it holds
\beam \label{eq.KP.5}
\limsup_{\xto} \dfrac{\bF(v\,x)}{\bF(x)} < 1\,.
\eeam

The following distribution class appears initially in \cite{konstantinides:tang:tsitsiashvili:2002}. We say that a distribution $F$ belongs to the distribution class $\mathcal{A}$, if $F \in \mathcal{S}$ and  \eqref{eq.KP.5} holds for any $v>1$, that means $\mathcal{A}:=\mathcal{S}\cap \mathcal{P_D}$. This class admits some useful properties in comparison with the subexponential one. Next we need the Matuszewska indexes, introduced in \cite{matuszewska:1964}. The  lower and upper Matuszewska indexes are given by
\beao
\beta_F:=\sup\left\{ -\dfrac{\ln \bF^*(v)}{\ln v} \;:\; v>1 \right\}\,,\qquad \alpha_F:=\inf\left\{ -\dfrac{\ln \bF_*(v)}{\ln v} \;:\; v>1 \right\}\,, 
\eeao
where
\beao
\bF_*(v):=\liminf_{\xto} \dfrac{\bF(v\,x)}{\bF(x)}\,, \qquad \bF^*(v):=\limsup_{\xto} \dfrac{\bF(v\,x)}{\bF(x)}\,,
\eeao
therefore, any distribution $F$, with right endpoint $r_F:=\sup\{y\;:\;F(y)<1\}=\infty$, satisfies the inequalities $0\leq \beta_F \leq \alpha_F$. It is known that the inclusion $F \in \mathcal{D}$ is equivalent to $\alpha_F < \infty$, and the inclusion $F \in \mathcal{P_D}$  is equivalent to $0<\beta_F$. Hence in class $\mathcal{A}$ the inequality $\beta_F>0$ is true, but it is not in $\mathcal{S}\setminus \mathcal{P_D}$. The restriction of the subexponentiality to positive decrease is frequent in the applications, see \cite{tang:2006}. For further results on class $\mathcal{A}$, see \cite{wang:hu:yang:wang:2018}, \cite{wang:zhang:wang:wang:2018}  and \cite{bardoutsos:konstantinides:2011}.

Now we consider a generalization of the class $\mathcal{L}$, introduced in \cite{shimura:watanabe:2005}. We say  $F$ belongs to the class of generalized long tailed distributions, symbolically $F \in \mathcal{OL}$, if for all (or equivalently, for some) $a > 0$ it holds
\beao
\limsup_{\xto}\dfrac{\bF(x-a)}{\bF(x)} < \infty\,.
\eeao

The generalization of the subexponential class was introduced in \cite{klueppelberg:1990}. We say that $F$ belongs to the class generalized subexponential distributions, symbolically $ F \in \mathcal{OS}$, if 
\beao
\limsup_{\xto}\dfrac{\overline{F^{2*}}(x)}{\bF(x)} < \infty\,.
\eeao
Further, in \cite{konstantinides:leipus:siaulys:2022}, the class $\mathcal{OA}:=\mathcal{OS}\cap \mathcal{P_D}$ was introduced. 

For a distribution $F$ with support $[0,\,\infty)$ and positive, finite expectation, $0<\E[X]< \infty$, we  consider the integrated tail distribution and its tail in the form
\beao
F_I(x):=\dfrac 1{\E[X]}\,\int_0^x \bF(y)\,dy\,,\qquad \bF_I(x)=\dfrac 1{\E[X]}\,\int_x^\infty \bF(y)\,dy\,,
\eeao
where $\E[X]$ denotes the mean value of $X$. In this case, from \cite{su:tang:2003} we find the classes $\mathcal{M}$ and $\mathcal{M}^*$. We say that $F$ belongs to the distribution class $\mathcal{M}$ if it holds
\beao
\lim_{\xto}\dfrac {\bF(x)}{\E[X]\,\bF_I(x)}=0\,.
\eeao
We say that $F$ belongs to distribution class $\mathcal{M}^*$ if holds
\beao
\limsup_{\xto}\dfrac {x\,\bF(x)}{\E[X]\,\bF_I(x)}< \infty\,.
\eeao
It is known that $\mathcal{M}^* \subseteq \mathcal{M} \subsetneq \mathcal{K}$. Furthermore, $\mathcal{D}\subsetneq \mathcal{M}^*$ (see \cite{su:tang:2003}). 

The following distribution class provides a multivariate version of regularly varying distributions and is called multivariate regular variation, symbolically $MRV$. We say that a random vector ${\bf X}$ with a multivariate distribution ${\bf F}$ belongs to class  $MRV$, if there exists a distribution $F \in \mathcal{R}_{-\alpha}$, with some $\alpha \in (0,\,\infty)$ and some Radon measure $\nu$ on the space $[0,\,\infty]^n\setminus\{(\infty,\,\ldots,\,\infty)\}$ such that
\beao
\lim_{\xto} \dfrac {\PP[{\bf X} \in x\,\bbb]}{\bF(x)} = \nu(\bbb)\,,
\eeao
as $\xto$, where $\bbb$ is any $\nu$-continuous Borel set from space $[0,\,\infty]^n\setminus\{{\bf 0}\}$. In this case we denote by ${\bf F} \in MRV(\alpha,\,F,\,\nu)$. The measure $\nu$ has the homogeneity property, means that for any $\bbb$ Borel set in space $[0,\,\infty]^n\setminus\{{\bf 0}\}$ and any $y>0$, the relation $\nu(y\,\bbb) = y^{-\alpha}\,\nu(\bbb)$ holds. This class  was introduced in \cite{dehaan:resnick:1982},  and triggered several applications in risk theory, see for example \cite{cheng:konstantinides:wang}, \cite{cheng:konstantinides:wang:2024}, \cite{konstantinides:li} and in risk management, see for example \cite{li:2022}, \cite{li:2022(b)}. Although class $MRV$ was well-established in many areas of applied probability,  the rest of multivariate distributions with heavy tails are still lesser known. 

In section 2, we provide some preliminaries known results in the case of independence and we extend these results  in the case of dependence, under some restrictions. 

In section 3 we introduce a new distribution class with some discussion on it, about the product distribution and related closure properties. Further, in section 4 we give some examples based on our main results. 

Additionally, in section 5 and 6 we introduce the dominatedly varying and positively decreasing random vectors respectively, and we study the closure property of the scalar product of each class under independent and dependent cases.

\section{Distribution of the product} \label{sec.KP.2}

In \cite{breiman:1965} and \cite{embrechts:goldie:1980} the  closure property of the product distribution in the class of regularly varying distribution tails was studied, providing some sufficiency conditions. In \cite{cline:samorodnitsky:1994} found some results for the distribution class $\mathcal{D}\cap \mathcal{L}$ with respect to the product distribution and some sufficient conditions for the inclusion $H \in \mathcal{S}$. However, the sufficient condition for the closure property is not mild enough for practical purposes. In \cite{tang:2006} a relaxation of these assumptions (see Assumption \ref{ass.KP.A}) and a reduction to class $\mathcal{A}$ facilitated a more useful result. Next, in \cite{xu:cheng:wang:chang:2018} we find a necessary and sufficient condition for $H \in \mathcal{S}$, when $F \in \mathcal{S}$, in the non-negative case, under some even more general condition than Assumption \ref{ass.KP.A} below. 

Next, we provide some results for the independent case, having in mind their further extension, under a specific weak dependence structure. Let us begin with the product convolution for the classes $\mathcal{M}$, $\mathcal{M}^*$ and $\mathcal{OL}$, as presented in \cite{su:chen:2007} and \cite{cui:wang:2020}.

\bpr [Su-Chen]\label{pr.Su-Chen.07}
Let $X$ and $Y$ be independent, non-negative random variables, with distributions $F$ and $G$ respectively, where $F$ has support $[0,\,\infty)$, and assume that $0<\E[Y]<\infty$. 
\begin{enumerate}

\item
If $F \in \mathcal{M}^*$, then it holds $H \in \mathcal{M}^*$. 

\item
If $F \in \mathcal{M}$, then it holds $H \in \mathcal{M}$. 
\end{enumerate}
\epr

\bpr [Cui-Wang]\label{pr.Cui-Wang.20}
Let $X$ and $Y$ be independent random  variables, with distributions $F$ and $G$ respectively. If $F \in \mathcal{OL}$, with support $[0,\,\infty)$, and  we assume $G(0-)=0$ and $G(0)<1$, then $H \in \mathcal{OL}$. 
\epr

In previous propositions  non-negative random variables only appear. We should remark that Proposition \ref{pr.Cui-Wang.20} appeared in \cite{cui:wang:2020}, in case when the support of $G$ is unbounded from above, and in \cite[Prop. 5.7(i)]{leipus:siaulys:konstantinides:2023}, in case when the support of $G$ is unbounded. In \cite{konstantinides:leipus:siaulys:2022} studied the distribution $F$ with support on the whole real axis, instead of class $\mathcal{OL}$ in the smaller class $\mathcal{OS}$, in order to achieve the closure property in $\mathcal{OS}$ under some additional condition. 

Before reminding \cite{konstantinides:leipus:siaulys:2022}, let us recall an assumption, that was used in several applications as sufficient condition for the closure property of some classes with respect to product.

\begin{assumption} \label{ass.KP.A}
For any $c>0$ it holds $\bG(c\,x)=o\left[\bH(x)\right]$, as $\xto$.
\end{assumption}

\bre \label{rem.KP.1}
From \cite[Lem. 3.2]{tang:2006} we find that when $G$ and $H$ satisfy $\bG(x)>0$ and $\bH(x)>0$ for any $x\geq 0$, then Assumption \ref{ass.KP.A} holds if and only if there exists a function $a\;:\;[0,\,\infty) \longrightarrow (0,\,\infty)$ such that
\begin{enumerate}

\item
$a(x) \longrightarrow \infty$,

\item
$a(x)/x \longrightarrow 0$,

\item
$\bG[a(x)] =o\left[\bH(x)\right]$,
\end{enumerate}
as $\xto$. It is remarkable that the $(1), (2), (3)$ represent the three of the four conditions in \cite[Th. 2.1]{cline:samorodnitsky:1994}. It is clear that, if the support of $F$ has infinite right endpoint, and the support of $G$ has finite right endpoint, then Assumption \ref{ass.KP.A} holds.
\ere

The following result can be found in \cite[Th. 3, Th. 4]{konstantinides:leipus:siaulys:2022}.

\bpr [Konstantinides-Leipus-\v{S}iaulys]\label{pr.K-L-S.22}
Let $X$ and $Y$ be independent random variables, with distributions $F$ and $G$ respectively, and let assume that $G(0-)=0\,,\; G(0)<1$, and $F$ has as support the whole real axis.
\begin{enumerate}

\item
If $F \in \mathcal{OS}$ and it holds
\beam \label{eq.KP.13}
\sup_{c>0}\,\limsup_{\xto} \dfrac{\bG(c\,x)}{\bH(x)} < \infty\,,
\eeam 
then it holds $H \in \mathcal{OS}$. 

\item
If $F \in \mathcal{OA}$ and Assumption \ref{ass.KP.A} is true, then it holds $H \in \mathcal{OA}$. 
\end{enumerate}
\epr

Later, \cite{mikutavicius:siaulys:2023} were able to show the closure property of class $\mathcal{OS}$ avoiding also relation \eqref{eq.KP.13}.

A plausible question comes up about the strictness of the independence condition between the random variables $X$ and $Y$. For example, when $X$ represents the net losses (claims minus premiums) and $Y$ represents the discount factor, it is reasonable to assume that if there is a 
large enough loss, this can affect the discount factor also. Next, we study the closure properties with respect to the product distribution on several heavy-tailed distribution classes, under a wide dependence frame. This dependence frame, introduced in \cite{asimit:jones:2008}, is given as follows.

\begin{assumption} \label{ass.KP.B}
Let exists a function $h\;:\;[0,\,\infty) \longrightarrow (0,\,\infty)$ such that the following relation $\PP[X>x\,|\,Y=y] \sim h(y)\,\PP[X > x]$, as $\xto$, holds uniformly for any $y \in S_Y$.
\end{assumption}

The uniformity is understood in the following sense:
\beao
\lim_{\xto}\,\sup_{y\in S_Y} \left| \dfrac{\PP[X>x\;|\;Y=y]}{h(y)\,\PP[X > x]} -1\right|=0\,.
\eeao
We can see that $\E[h(Y)]=1$ holds. Further, we observe that function $h(y)$ is bounded from above for any $y\geq 0$ (see \cite[Lem. 3.1]{chen:xu:cheng:2019}).

This dependence frame includes a variety of copulas, as for example the Farlie-Gumbel-Morgenstern copula or the Frank copula. This dependence belongs to the family of asymptotic independence, and in the discrete time risk models describes the dependence between main loss and discount factor, while in continuous time risk model the dependence between the claim sizes and the inter-arrival times (see   \cite{yang:gao:li:2016}, \cite{yang:leipus:siaulys:2012} and \cite{asimit:badescu:2010}, \cite{li:tang:wu:2010} respectively). In general there are several papers about the closure properties of the product distributions for heavy-tailed distributions. In \cite{yang:wang:2013} we find the closure properties of some distribution classes under Sarmanov dependence. In \cite{cadena:omey:vesilo:2022}, \cite{cui:wang:2024} they studied the closedness under the product for some heavy-tailed distribution classes, under Assumption \ref{ass.KP.B}. There are some papers for random vectors from heavy-tailed distributions, where this problem is reduced to study of the closure properties of random vector product distribution (see \cite{basrak:davis:mikosch:2002}, \cite{samorodnitsky:sun:2016} for independent factors and \cite{fougeres:mercadier:2012} for dependent factors).

Now we begin with \cite[Lem. 3.1]{chen:xu:cheng:2019}, needed for our further results.

\ble \label{lem.Chen-Xu-Cheng.19}
Let $X$ and $Y$ be non-degenerated random variables with distributions $F$ and $G$ respectively. Let the support of $F$ be the whole real axis and for $G$ let us assume that $G(0-)=0$ and $G(0)<1$. If $X$ and $Y$ satisfy the Assumptions \ref{ass.KP.A} and \ref{ass.KP.B}, then it holds
\beao
\PP[X\,Y > x] \sim \int_{0}^\infty h(y)\,\bF\left(x/y\right)\,G(dy)\,,
\eeao
as $\xto$.
\ele

The next result represents an extension, in case of weak dependence structure of Assumption \ref{ass.KP.B}, of Proposition \ref{pr.Su-Chen.07}. The extra restriction we use, in comparison with the independent case, is Assumption \ref{ass.KP.A}.

\bth \label{th.KP.2}
Let $X$ and $Y$ be non-negative random variables with distributions $F$ and $G$ respectively, and $F$ has as support $[0,\,\infty)$. If Assumptions \ref{ass.KP.A} and \ref{ass.KP.B} hold and $0<\E[Y]<\infty$, we obtain
\begin{enumerate}

\item
If $F \in \mathcal{M}^*$, then it holds $H \in \mathcal{M}^*$. 

\item
If $F \in \mathcal{M}$, then it holds $H \in \mathcal{M}$. 
\end{enumerate}
\ethe

\pr~
\begin{enumerate}

\item

It is proved that there exists $K>0$ such that
\beam \label{eq.KP.22}
h(y) \leq K\,,
\eeam
for any $y\geq 0$, see in \cite[Th. 2.1]{chen:xu:cheng:2019} or \cite[Prop. 2.4]{cui:wang:2024}. Now we define a random variable $Y_h$, that is independent of $X$, and has distribution
\beam \label{eq.KP.23}
G_h(dy) =h(y)\,G(dy)\,,
\eeam
for any $y\geq 0$, where the $G_h$ is proper distribution, because of $\E[h(Y)]=1$. Hence, the random variables $X$ and $Y_h$ satisfy the conditions of Proposition  \ref{pr.Su-Chen.07}. Indeed, by \eqref{eq.KP.22} we obtain
\beao
\E[Y_h]=\int_0^\infty y\,G_h(dy)=\int_0^\infty y\,h(y)\,G(dy)\leq K\,\int_0^\infty y\,G(dy) \leq K\,E[Y]< \infty\,,
\eeao
and $\E[Y_h]>0$, since function $h$ is strictly positive by definition. Finally, applying Lemma \ref{lem.Chen-Xu-Cheng.19} we find
\beam \label{eq.KP.24}
\bH(x) := \PP[X\,Y>x] \sim \PP[X\,Y_h > x]=:\bH_h(x)\,,
\eeam
as $\xto$, whence by applying Proposition \ref{pr.Su-Chen.07}(1) on $X$ and $Y_h$, we obtain $H_h \in \mathcal{M}^*$, which together with \eqref{eq.KP.24} provides $H \in \mathcal{M}^*$.

\item
In similar way, we apply Proposition \ref{pr.Su-Chen.07}(2)  on $X$ and $Y_h$.~\halmos
\end{enumerate}

In the next result we take into account the class of positive decreasing distributions.

\bth \label{th.KP.3}
Let $X$ and $Y$ be random variables with distributions $F$ and $G$ respectively. The distribution $F$ has as support the whole real axis and the distribution $G$ satisfies the conditions $G(0-)=0\,,\; G(0)<1$. Under the Assumptions \ref{ass.KP.A} and \ref{ass.KP.B}, if $F \in \mathcal{P_D}$, then it holds $H \in \mathcal{P_D}$.  
\ethe

\pr~
Let $\bG(x)>0$, for any $x\geq 0$. For any $v>1$, applying Lemma \ref{lem.Chen-Xu-Cheng.19} we get
\beao
\limsup_{\xto}\dfrac{\bH(v\,x)}{\bH(x)}&=&\limsup_{\xto}\dfrac 1{\bH(x)}\left(\int_0^{a(x)} + \int_{a(x)}^\infty \right) \PP\left[X >\dfrac{v\,x}y\;|\;Y=y \right]G(dy)\\[2mm]
&\leq& \limsup_{\xto} \dfrac{\int_0^{a(x)}\PP\left[X>\dfrac {v\,x}y\;|\;Y=y\right]\,G(dy)}{\int_0^{a(x)}\PP\left[X>\dfrac xy\;|\;Y=y \right]\,G(dy)} +\limsup_{\xto} \dfrac{\bG[a(x)]}{\bH(x)}\\[2mm]
&=&\limsup_{\xto} \int_0^{a(x)}\dfrac{ h(y)\,\bF\left(v\,x/y\right)}{\int_0^{a(x)}h(y)\,\bF\left(x/y \right)\,G(dy)}\,G(dy) \\[2mm]
&\leq& \limsup_{\xto} \sup_{0< y \leq a(x)} \dfrac{\bF\left({v\,x}/y\right)}{\bF\left(x/y \right)}= \limsup_{\xto} \dfrac{\bF\left(v\,x \right)}{\bF\left( x \right)} < 1\,,
\eeao
where $a(x)$ follows by Assumption \ref{ass.KP.A} (see Remark \ref{rem.KP.1}). In the second step used the condition $\bG[a(x)]=o[\bH(x)]$, as $\xto$, and in the last step we take into account that $F \in \mathcal{P_D}$. Therefore, we obtain $H \in \mathcal{P_D}$. In case the support of $G$ is bounded from above, the proof follows similarly, putting instead of $a(x)$, the right endpoint of the support of $G$.
~\halmos

In the next result we follow Proposition \ref{pr.K-L-S.22} under the dependence of Assumption \ref{ass.KP.B}.

\bth \label{th.KP.4}
Let $X$ and $Y$ be random variables with distributions $F$ and $G$ respectively. The support of $F$ is the whole real axis and for $G$ we assume $G(0-)=0\,,\; G(0)<1$. Let assume that Assumptions \ref{ass.KP.A} and \ref{ass.KP.B} are true and
\begin{enumerate}

\item
if $F \in \mathcal{OS}$, then it holds $H \in \mathcal{OS}$. 

\item
if $F \in \mathcal{OA}$, then it holds $H \in \mathcal{OA}$. 
\end{enumerate}
\ethe

\pr~
\begin{enumerate}

\item
We consider a random variable $Y_h$ with distribution given in \eqref{eq.KP.23},  and independent of $X$. Hence, by Assumptions \ref{ass.KP.A} and \ref{ass.KP.B} and Lemma   \ref{lem.Chen-Xu-Cheng.19}, we obtain relation \eqref{eq.KP.24}. Further, we observe that Assumption \ref{ass.KP.A} is included in  \eqref{eq.KP.13}, and since for any $c>0$ it holds
\beao
\bG_h(c\,x)=\int_{c\,x}^\infty G_h(dy)=\int_{c\,x}^\infty h(y)\,G(dy) \leq K\,\bG(c\,x)\,,
\eeao
for some constant $K>0$, then it follows
\beam \label{eq.KP.2.5}
\bG_h(c\,x)=o\left[ \bH(x) \right]=o\left[ \bH_h(x) \right]\,,
\eeam
for any constant $c>0$, where in the last step was used \eqref{eq.KP.24}. Therefore, by Proposition \ref{pr.K-L-S.22}(1) we obtain $H_h \in \mathcal{OS}$. And further by \eqref{eq.KP.24} we conclude $H \in \mathcal{OS}$.

\item
From the first part we obtain $H \in \mathcal{OS}$. Further, by Theorem \ref{th.KP.3} is implied that $H \in \mathcal{P_D}$, thus we have  $H \in \mathcal{OA}$.~\halmos
\end{enumerate}

\bre \label{rem.KP.2}
In Theorem \ref{th.KP.4} is clearly demonstrated that in class $\mathcal{OS}$ we need a
stricter condition for the case of Assumption \ref{ass.KP.B},  in comparison to the independence (note that Assumption \ref{ass.KP.A} is contained in relation \eqref{eq.KP.13}), while in class $\mathcal{OA}$ the conditions remain the same. Furthermore,  the situation can NOT change in case of reduction of the support of $F$ in the interval $[0,\,\infty)$, since Assumption \ref{ass.KP.A} is necessary for the application of Lemma  \ref{lem.Chen-Xu-Cheng.19}.
\ere

In the following result we consider the distribution class $\mathcal{OL}$, only when $F$ has as support the interval $[0,\,\infty)$. As we have seen in Proposition \ref{pr.Cui-Wang.20}, do NOT exist conditions for inclusion $H \in \mathcal{OL}$, when $F \in \mathcal{OL}$, in case of independence. However in the next result we have the dependence under the framework of Assumption  \ref{ass.KP.B}, whence we have to use Assumption \ref{ass.KP.A}.

\bth \label{th.KP.5}
Let $X$ and $Y$ be  non-negative  random variables with distributions $F$ and $G$ respectively, where we assume that $G(0-)=0\,,\; G(0)<1$. If Assumptions \ref{ass.KP.A} and \ref{ass.KP.B} are true and $F \in \mathcal{OL}$ then it holds $H \in \mathcal{OL}$.  
\ethe 

\pr~
We consider the independent random variables $X$ and $Y_h$, with the distribution of $Y_h$ given by \eqref{eq.KP.23}. From Lemma  \ref{lem.Chen-Xu-Cheng.19} we obtain that \eqref{eq.KP.24} is true, thence by Proposition \ref{pr.Cui-Wang.20} we find $H_h \in \mathcal{OL}$ and therefore $H \in \mathcal{OL}$.
~\halmos

The following result gives closure properties of the product distribution for the classes $\mathcal{D}$, $\mathcal{D}\cap \mathcal{L}$, $\mathcal{C}$ and $\mathcal{K}$ under the dependence of Assumption \ref{ass.KP.B}, and under Assumption \ref{ass.KP.A}. In \cite{yang:wang:leipus:siaulys:2013} we find a similar result for the distribution class $\mathcal{D}$ under the condition 
\beam \label{eq.KP.25}
\bG(x)= o\left[\bF(x)\right]\,,
\eeam
as $\xto$, where Assumption \ref{ass.KP.A} is included in condition \eqref{eq.KP.25}, see Remark \ref{rem.KP.1}. The second part in the next statement corresponds to \cite[Th. 2.2(iii)]{cline:samorodnitsky:1994}, in case of dependence in the frame of Assumption \ref{ass.KP.B}. It is worth to notice that there is no extra condition with respect to distributions, in comparison with the independent case. Further, the second part covers \cite[Lem 2.1]{yang:leipus:siaulys:2012} and \cite[Th. 3.3]{cadena:omey:vesilo:2022}, since  Assumption \ref{ass.KP.A} is more general. Even more, parts 4 and 5 correspond to \cite[Th. 3.4(ii) and Th. 3.5(ii)]{cline:samorodnitsky:1994}, in case of dependence as in the frame of Assumption \ref{ass.KP.B}.

\bth \label{th.KP.6}
Let $X$ and $Y$ be random variables with distributions $F$ and $G$ respectively, where $F$ has as support the whole real axis and for $G$ we assume $G(0-)=0\,,\; G(0)<1$. We assume that Assumptions \ref{ass.KP.A} and \ref{ass.KP.B} are true and
\begin{enumerate}
\item
if $F \in \mathcal{D}$, then it holds $H \in \mathcal{D}$.
\item
if $F \in \mathcal{D}\cap \mathcal{L}$, then it holds $H \in \mathcal{D}\cap \mathcal{L}$. 
\item 
if $F\in \mathcal{C}$, then it holds $H \in \mathcal{C}$. 
\item 
if $F \in\mathcal{K}$, then it holds $H \in\mathcal{K}$.
\end{enumerate}
\ethe

\pr~
\begin{enumerate}
\item
As in Theorem \ref{th.KP.3}, we show only the case $\bG(x)>0$, for any $x\geq 0$, since the proof in case the support of $G$ has finite right endpoint follows similarly by putting instead of $a(x)$ this right endpoint of the support. Let $0<b<1$, then
\beao
\limsup_{\xto} \dfrac{\bH(b\,x)}{\bH(x)}&=&\limsup_{\xto} \dfrac 1{\bH(x)}\left(\int_0^{a(x)}+\int_{a(x)}^\infty \right) \PP\left[X>\dfrac {b\,x}y\;\Big|\;Y=y \right]\,G(dy) \\[2mm]
&\leq& \limsup_{\xto} \dfrac{ \int_0^{a(x)}\PP\left[ X>\dfrac{b\,x}y\;\Big|\;Y=y\right]\,G(dy)}{\int_0^{a(x)} \PP\left[X>\dfrac {x}y\;\Big|\;Y=y \right]\,G(dy)} + \limsup_{\xto} \dfrac{\bG[a(x)]}{\bH(x)} \\[2mm]
&=& \limsup_{\xto} \int_0^{a(x)}\dfrac{ h(y)\,\bF\left(b\,x/y\right)}{\int_0^{a(x)} h(y)\,\bF\left(x/y\right)\,G(dy)}\,G(dy)  \\[2mm]
&\leq& \limsup_{\xto} \,\sup_{0< y \leq a(x)}\dfrac{\bF\left(b\,x/y\right)}{\bF\left( x/y \right)}= \limsup_{\xto} \dfrac{\bF\left(b\,x\right)}{\bF\left(x \right)} < \infty\,,
\eeao 
where $a(x)$ denotes the function in Remark \ref{rem.KP.1}, we used Lemma \ref{lem.Chen-Xu-Cheng.19} and in the last step we take into account $F \in \mathcal{D}$. Therefore, we obtain $H \in \mathcal{D}$.

\item
This follows from previous part (1) and from \cite[Th. 2.1(i)]{chen:xu:cheng:2019}.

\item 
Consider that $F \in \mathcal{C}$ and the random variable $Y_h$ is independent of $X$, with distribution given by \eqref{eq.KP.23}. By  Assumptions \ref{ass.KP.A} and \ref{ass.KP.B},  we obtain  \eqref{eq.KP.24} and next we get relation \eqref{eq.KP.2.5}. From \cite[Th. 3.4(ii)]{cline:samorodnitsky:1994} we find that $H_h \in \mathcal{C}$ and thus we have that $H \in \mathcal{C}$.

\item
Now we consider that $F \in \mathcal{K}$ and the random variable $Y_h$ is independent of $X$, with distribution given by \eqref{eq.KP.23}. By of Assumptions \ref{ass.KP.A} and \ref{ass.KP.B}, we obtain  \eqref{eq.KP.24}. From \cite[Prop. 5.1(i)]{leipus:siaulys:konstantinides:2023} we find that $H_h \in \mathcal{K}$ and thus we have $H \in \mathcal{K}$.
~\halmos
\end{enumerate}

The following result studies the closure property of product distribution under independent and dependent cases of class $\mathcal{D}\cap\mathcal{A}$. It is easy to see that $\mathcal{D}\cap\mathcal{A}:=\mathcal{D}\cap\mathcal{S}\cap\mathcal{P_{D}}=\mathcal{D}\cap\mathcal{L}\cap\mathcal{P_{D}}$, which follows from the fact that $\mathcal{D}\cap\mathcal{S}=\mathcal{D}\cap\mathcal{L}$, see \cite[Th.1]{goldie:1978}.

\bco\label{cor.KP.1}
Let $X$ and $Y$ be random variables with distributions $F$ and $G$ respectively. Distribution $F$ has support the whole real axis and for $G$ we assume $G(0-)=0\,,\; G(0)<1$. If Assumptions \ref{ass.KP.A} and \ref{ass.KP.B} hold, and $F \in \mathcal{D}\cap\mathcal{A}$ is true, then it holds $H \in \mathcal{D}\cap\mathcal{A}$. 
\eco

\pr~
By \cite[Th.3.1]{wang:zhang:wang:wang:2018} we get that $ H \in\mathcal{A}$, so by Theorem \ref{th.KP.6}(1) we have that $H \in \mathcal{D}$.~\halmos

\section{Distribution class $\mathcal{T}$} \label{sec.KP.4}

In this section we introduced a new class of heavy-tailed distributions. This class have some remarkable properties, as for example the closure properties (i) with respect to the distribution product (in both independent and dependent cases), (ii) with respect to convolution of independent random variables, (iii) with respect to the mixture of arbitrarily dependent random variables.

\bde \label{def.KP.1}
We say that a distribution $F$ with support the whole real axis, belongs to the long tailed positively decreasing distributions, symbolically $\mathcal{T}$, if $F \in \mathcal{L}$, and $F \in \mathcal{P_D}$, symbolically $\mathcal{T}:=\mathcal{L}\cap\mathcal{P_D}$.
\ede

It is obvious that, if $F \in\mathcal{T}$, then it follows $\beta_F>0$, $\mathcal{A}\subset\mathcal{T}$, and $\mathcal{D}\cap\mathcal{A}=\mathcal{D}\cap\mathcal{T}$.

\bth \label{th.KP.7}
Let $X$ and $Y$ be random variables, with distributions $F$ and $G$ respectively, $F$ has support of the whole real axis and for $G$ we assume 
$G(0-)=0\,,\; G(0)<1$. If $X$ and $Y$ satisfy Assumption \ref{ass.KP.A} and \ref{ass.KP.B} and $F \in \mathcal{T}$ is true, then it holds $H \in \mathcal{T}$. 
\ethe

\pr~
By combination of Theorem \ref{th.KP.3} and \cite[Th. 2.1(i)]{chen:xu:cheng:2019}, we obtain the desired result.~\halmos 

The following result renders closure properties of class $\mathcal{T}$, which are not relevant to the concept of the product distribution. For this reason we need some extra notations. Let remind for some random variables $X_{1},\,X_{2}$, with distributions $F_1,\,F_2$ respectively, which belong to class $\mathcal{B}$, and now considering the  $X_{1},\,X_{2}$ are arbitrarily dependent, we say that this class satisfies the closure property with respect to mixture if $p\,F_{1}+(1-p)\,F_{2} \in \mathcal{B}$, holds for any $p \in (0,\,1)$. 

\bpr \label{th.KP.8}
\begin{enumerate}
\item 
Let  $X_{1},\,X_{2}$ be arbitrarily dependent random variables with distributions $F_{1},\,F_{2} \in \mathcal{P_{D}}$ respectively, then, for any $p \in (0,\,1)$, it holds $p\,F_{1}+(1-p)\,F_{2} \in \mathcal{P_{D}}$.  

\item 
Let  $X_{1},\,X_{2}$ be arbitrarily dependent random variables with $F_1 \in\mathcal{T}$ and either $F_2 \in\mathcal{T}$, or $\bF_2 (x)=o\left[\overline{F}_{1}(x)\right]$ and $F_2 \in\mathcal{P_{D}}$, then, for any $p \in (0,\,1)$, it holds $p\,F_{1}+(1-p)\,F_{2} \in \mathcal{T}$. 
\end{enumerate}
\epr

\pr~
\begin{enumerate}
\item 
Let take arbitrarily chosen $v>1$. Then for any $p \in (0,\,1)$ we obtain 
\beao
	\limsup_{\xto}\dfrac{p\,\overline{F}_{1}(v\,x)+(1-p)\,\overline{F}_{2}(v\,x)}{p\,\overline{F}_{1}(x)+(1-p)\,\overline{F}_{2}(x)} \leq \max\left\{ \limsup_{\xto}\dfrac{\bF_{1}(v\,x)}{\bF_{1}(x)}, \dfrac{\bF_{2}(v\,x)}{\bF_{2}(x)} \right\}<1\,,
\eeao
hence $p\,F_{1}+(1-p)\,F_{2} \in \mathcal{P_{D}}$.

\item 
In the case $F_{1} \in \mathcal{T}$ and $F_{2} \in \mathcal{T}$, by \cite[Cor. 2.23]{Foss:Korshunov:Zachary:2013} we have that since $F_{1} \in \mathcal{L}$ and $F_{2} \in \mathcal{L}$, then it holds $p\,F_{1}+(1-p)\,F_{2} \in\mathcal{L}$. This result, together with part (1), implies $p\,F_{1}+(1-p)\,F_{2} \in \mathcal{T}$.

In the other case, where $F_{1}\in\mathcal{T}$ and $\overline{F}_{2}(x)=o\left[\overline{F}_{1}(x)\right]$, with $F_{2} \in \mathcal{P_D}$, by \cite[Prop. 3.9]{leipus:siaulys:konstantinides:2023} we obtain that $p\,F_{1}+(1-p)\,F_{2} \in \mathcal{L}$. Their combination, together with part (1), renders the required result.
~\halmos
\end{enumerate}

\section{Applications}

Next, we provide some applications that demonstrate the previous results on randomly weighted sums and ruin probabilities in a discrete-time risk model with actuarial and financial risks. Let us introduce the notation $S_n^Y := \sum_{i=1}^n Y_i\,X_i\,, \; M_n^Y := \max_{1\leq i \leq n} S_i^Y$. The last case corresponds to the ruin probability in the risk model over discrete time, with stochastic discount factors. Let us notice that in present section, when we say that Assumption \ref{ass.KP.B} is satisfied for any $(X_i,\,Y_j)$, then Assumption  \ref{ass.KP.B} can be satisfied with $h_i \neq h_j$ for $1\leq i \neq j \leq n$.

In the first example we obtain a generalization of \cite[Th. 1]{yang:wang:leipus:siaulys:2013}, since the Assumption \ref{ass.KP.A} is more general in comparison with their assumption, see \cite[Cor. 2.1]{tang:2006}.

\bexam \label{exam.KP.1}
Let $(X_1,\,Y_1),\,\ldots,\,(X_n,\,Y_n)$ be mutually independent random vectors, where the $X_1,\,\ldots,\,X_n$ are real random variables with distributions $F_1,\,\ldots,\, F_n$ respectively and the $Y_1,\,\ldots,\,Y_n$ strictly positive random variables with distributions $G_1,\,\ldots,\,G_n$ respectively. Let assume that the pair $(X_i,\,Y_i)$, satisfies Assumptions  \ref{ass.KP.A} and \ref{ass.KP.B}, for any $i=1,\,\ldots,\,n$.
\begin{enumerate}

\item
If  $F_i \in \mathcal{L}$, for  $i=1,\,\ldots,\,n$, then it holds
\beao
\PP\left[S_n^Y >x\right] \sim \PP\left[M_n^Y >x\right] \sim \PP\left[\sum_{i=1}^n Y_i\,X_i^+ >x\right]\,,
\eeao
as $\xto$.

\item
If $F_i \in \mathcal{D}\cap \mathcal{L}$, for  $i=1,\,\ldots,\,n$, then it holds
\beao
\PP\left[S_n^Y >x\right] \sim \PP\left[M_n^Y >x \right] \sim \PP\left[\sum_{i=1}^n Y_i\,X_i^+ >x\right] \sim \sum_{i=1}^n \PP\left[Y_i\,X_i >x\right]\,,
\eeao
as $\xto$.
\end{enumerate}
\eexam

\pr~
We employ the argument from \cite[Th. 1]{yang:wang:leipus:siaulys:2013}, but instead of their Lemmas 1 and 2, we use \cite[Th. 2.1(i)]{chen:xu:cheng:2019} for $\gamma =0$ and Theorem \ref{th.KP.6}(1) respectively.
~\halmos

Next, we study the asymptotic behavior of the ruin probability, in a discrete-time risk model, in which the financial and actuarial risks satisfy the dependence of Assumption \ref{ass.KP.B}. If we denote by $X_i$ the net loss of the insurer during the $i$-th period, namely the claims minus premiums in this period, and by $Y_i$ the stochastic discount factor, which takes non-negative values, then the total discount value of the insurer's net losses during the first $n$ periods, is given by
\beam \label{eq.KP.4.A}  
S_n^{\Theta} = \sum_{i=1}^n X_i\,\Theta_i = \sum_{i=1}^n X_i\,\prod_{j=1}^i Y_j\,,
\eeam
for any $n \in \bbn \cup \{\infty\}$, with $S_0=0$ and $\Theta_i :=\prod_{j=1}^i Y_j$. Hence, if we consider that the insurer's initial capital is $x>0$, then the ruin probability over finite or infinite time horizon, in model \eqref{eq.KP.4.A} is defined by
\beam \label{eq.KP.4.B}
\psi(x,\,n):= \PP[M_n >x]\,,
\eeam
for any $n \in \bbn \cup \{\infty\}$, where
\beam \label{eq.KP.4.C}
M_n:=\max_{0\leq k \leq n} S_k^{\Theta}\,,
\eeam
for any  $n \in \bbn_0 \cup \{\infty\}$. The risk model in \eqref{eq.KP.4.A}, was examined first time in \cite{tang:tsitsiashvili:2003}, with independent $X_i$, $Y_i$. Since then appeared several dependent, risk models, see for example \cite{chen:2011}, \cite{yang:konstantinides:2015}, \cite{yang:gao:li:2016}, \cite{chen:yuan:2017}, \cite{tang:yang:2019} among others. In  \cite{chen:yuan:2017}, was studied the ruin probability of model in \eqref{eq.KP.4.A} with finite time horizon, under the condition that the product distribtuion $X\,Y$, let call it $H$, belongs to class $\mathcal{C}$, or to class $\mathcal{P_D}$ for infinite time horizon, under arbitrary dependence between  $X$ and $Y$, and under some moment conditions for $Y$, with respect to Matuszewska indexes of  $H$. 

In next result, we show that under Assumptions \ref{ass.KP.A} and \ref{ass.KP.B}, with $F$ from class $\mathcal{B} \in \{\mathcal{C},\,\mathcal{P_D}\}$, it follows that $H \in \mathcal{B}$, and further we relax the moment conditions on $Y$, from the Matuszewska indexes of  $H$, to Matuszewska indexes of  $F$.

\bco
\begin{enumerate}
\item
Let us consider the risk model in \eqref{eq.KP.4.A}. If the sequence $\{(X_i,\,Y_i)\,,\;i \in \bbn\}$ of independent, identically distributed copies of the couple $(X,\,Y)$, with distributions $F,\,G$, such that $F\in \mathcal{C}$,  $G(0-)=0$, $G(0)<1$ and further the couple  $(X,\,Y)$ satisfies  Assumptions \ref{ass.KP.A} and \ref{ass.KP.B}, and $\E[Y^{\alpha_F +\vep}] < \infty$ for some $\vep >0$, then it holds
\beam \label{eq.KP.4.D}
\psi(x,\,n) \sim  \sum_{i=1}^n\PP\left[ X_i\,\Theta_i >x\right]\,,
\eeam
as $\xto$.
\item
Under the assumption of the part (1), with the restriction that $F \in \mathcal{C}\cap \mathcal{P_D}$, and $\E[Y^{\alpha}\vee Y^{\beta}] < 1$, for some $0\leq \beta \leq \beta_F \leq \alpha_F \leq \alpha < \infty$, then \eqref{eq.KP.4.D} holds uniformly, for any $n \in \bbn$, and hence also for $n=\infty$. Namely,
\beam \label{eq.KP.4.E}
\lim_{\xto} \sup_{n \in \bbn} \left|\dfrac{\psi(x,\,n)}{ \sum_{i=1}^n\PP\left[ X_i\,\Theta_i >x\right]} -1 \right| =0\,.
\eeam
\end{enumerate}
\eco

\pr~
At first, if $F \in \mathcal{B}  \in \{\mathcal{C},\,\mathcal{P_D}\} $, then $H \in \mathcal{B}$, which follows by Theorem \ref{th.KP.6}(3) and Theorem \ref{th.KP.3}. Because of $\mathcal{B} \subsetneq \mathcal{D}$, and from the condition $\E[Y^{\alpha_F +\vep}] < \infty$, we obtain that from \cite[Lem. 3.9]{tang:tsitsiashvili:2003b}, for the independent random variable $Y_h$, whose distribution is $G_h$, see relation \eqref{eq.KP.23}, that satisfies the following relation
\beam \label{eq.KP.4.ST}
\E[Y_h^{\alpha_F +\vep}] =\int_{0}^{\infty} y^{\alpha_F + \vep}\,G_h(dy) \leq k\, \E[Y^{\alpha_F +\vep}]  <\infty\,,
\eeam
for some $\vep >0$, is implied the equalities $\alpha_F = \alpha_{H_n}$, $\beta_F=\beta_{H_n}$, the same holds and for the stricter moment conditions in part (2). Hence, since relation \eqref{eq.KP.24} holds, we find that $\alpha_H = \alpha_{H_n}$, $\beta_H=\beta_{H_n}$, thus we conclude in both cases $\alpha_F = \alpha_{H}$, $\beta_F=\beta_{H}$. Therefore, 
\begin{enumerate}
\item
we obtain $H \in \mathcal{C}$, with $\E[Y^{\alpha_H +\vep}] < \infty$, and applying \cite[Th. 3.2]{chen:yuan:2017}, we find the desired result.
\item
we obtain $H \in \mathcal{C}\cap \mathcal{P_D}$, and $\E[Y^{\alpha}\vee Y^{\beta}] < 1$, for some $0\leq \beta \leq \beta_F \leq \alpha_F \leq \alpha < \infty$, hence applying \cite[Th. 3.3]{chen:yuan:2017}, we find the desired result. 
\end{enumerate}
~\halmos

\section{Dominatedly varying random vectors} \label{sec.KP.5}

Recently the competition among the insurance companies and financial institutions created the need to study the heavy tailed distributions, not individually, but in vector form, where the dependence among components is permitted. The most known random vector of this kind, represented by the multivariate regular variation, was introduced in \cite{dehaan:resnick:1982}. In this case, the approximation  describes a random vector with regularly varying components, with identical regular variation index $\alpha \in (0,\,\infty)$ (see in \cite{resnick:2007} and in relation with risk theory applications in \cite{konstantinides:li}, \cite{li:2016}, \cite{cheng:konstantinides:wang}). 

There were several attempts to define random vectors in larger classes of heavy tailed distributions, mainly in multivariate subexponential distributions. In \cite{cline:resnick:1992} was proposed the multivariate subexponentiality through point processes. Later, in \cite{omey:2006} were given three different definitions of  multivariate subexponential distributions, in more direct way. Finally in \cite{samorodnitsky:sun:2016} was introduced another definition of multivariate subexponential distribution, which was accompanied by a continuous time risk model application. Our definitions in section 5 and 6  for the classes $\mathcal{D}$ and $\mathcal{P_D}$ are inspired by these from \cite{omey:2006} for the subexponential class $\mathcal{S}$, but there exists a basic difference. While in \cite{omey:2006} was considered the distribution tail of the form ${\bf \overline{F}}({\bf t}\,x):=1-{\bf F}({\bf t}\,x)=1-\PP[X_1 \leq t_1\,x,\,\ldots,\,X_n \leq t_n\,x]$, namely, at least on excess for $X_1,\,\ldots,\,X_n$, we suggest the distribution tail \eqref{eq.KP.5.1}, where  the excesses of the components happen all together.

Let consider a random vector ${\bf X}=(X_1,\,\ldots,\,X_n)$, whose distribution tail is defined by the formula
\beam \label{eq.KP.5.1}
{\bf \bF}({\bf t},\,x)=\PP[X_1>t_1\,x,\,\ldots,\,X_n>t_n\,x]\,,
\eeam
for any ${\bf t}=(t_1,\,\ldots,\,t_n) \in (0,\,\infty]^n \setminus \{(\infty,\,\ldots,\,\infty)\}$. We notice that, if for some component of ${\bf t}$, we have $t_i=\infty$, then we get automatically reduction of the dimension.

Let us introduce a new class of distributions and demonstrate its usage. The ${\bf t}$ plays the role of tail convergence and in fact helps the convergence flexibility, and its domain of definition was chosen similarly to that in \cite{omey:2006}. In relation to \cite{omey:2006}, which focus mainly on multivariate subexponential and multivariate long-tailed distributions, here the characterization of the classes is carried out only through the joint tail, having in mind that the immediateness with uni-variate definitions, and the presence of single big jump principle of each component as well. 

\bde \label{def.KP.5.1}
Let random vector ${\bf X}$ with distribution ${\bf F}$ and its marginal distributions $F_1,\,\ldots,\,F_n$. We assume that there exists a distribution $F$, such that $F \in \mathcal{D}$ and $\bF_i(x) \asymp \bF(x)$, as $\xto$, for any $i=1,\,\ldots,\,n$. If further, it holds
\beam \label{eq.KP.5.2} 
\limsup_{\xto}\dfrac{{\bf \bF}({\bf b}\,{\bf t},\,x)}{{\bf \bF}({\bf t},\,x)}=\limsup_{\xto}\dfrac{\PP[X_1>b_1\,t_1\,x,\,\ldots,\,X_n>b_n\,t_n\,x]}{\PP[X_1>t_1\,x,\,\ldots,\,X_n>t_n\,x]} <\infty\,,
\eeam
for all ${\bf t} \in (0,\,\infty]^n \setminus \{(\infty,\,\ldots,\,\infty)\}$ and all (or equivalently, for some) ${\bf b} \in (0,\,1)^n$, then we say that random vector ${\bf X}$ follows multivariate dominatedly varying distribution with weak kernel $F$ and we write ${\bf F} \in \mathcal{D}_n$.
\ede

\bre \label{rem.KP.5.1}
From condition $\bF_i(x) \asymp \bF(x)$, as $\xto$, for any $i=1,\,\ldots,\,n$, we obtain $F_i \in \mathcal{D}$, for $i=1,\,\ldots,\,n$, (see, for example \cite[Prop. 3.7(i)]{leipus:siaulys:konstantinides:2023}. Specifically, in class $\mathcal{D}_n$, there exist by definition a weak kernel.

In the case, when vector  ${\bf X}$ has non-negative components, a nice feature, which appears in multivariate dominatedly varying distributions, is the stability under non-negative and non-degenerated linear combinations. As mentioned in \cite{samorodnitsky:sun:2016}, many well established multivariate extensions of one dimensional definitions in probability theory, have the following characteristic: if random vector ${\bf X}=(X_{1},\,\ldots,\,X_{n})$ has the property $\mathcal{P}_{n}$ (which symbolizes the n-dimensional property), then any non-negative, nondegenerated to zero, linear combination $\sum_{i=1}^{n}l_{i}X_{i}$ has the property $\mathcal{P}_{1}$. For example the multivariate normal and the multivariate regular variation have this feature. Actually, in the non-negative case, the class of multivariate dominatedly varying distributions is equipped this device. Let  ${\bf X}=(X_{1},\,\ldots,\,X_{n})$ with distribution  ${\bf F} \in \mathcal{D}_n$, then for any linear combination and any $b\in(0,\,1)$ holds
\beao
	\limsup_{\xto}\dfrac{P\left[\sum_{i=1}^{n}l_{i}X_{i}>b\,x \right]}{P\left[\sum_{i=1}^{n}l_{i}X_{i}>x\right]}\leq\limsup_{\xto}\dfrac{P\left[\sum_{i=1}^{n}X_{i}>b\,\dfrac x{\check{l}}\right]}{P\left[\sum_{i=1}^{n}X_{i}> \dfrac x{\widehat{l}}\right]}<\infty\,,
\eeao
where $\check{l}:=\max\{l_{1},\,\ldots,\,l_{n}\}$, $\widehat{l}:=\min\{l_{1},\,\ldots,\,l_{n}\}$, and the last step follows by the closure property of class $\mathcal{D}$ with respect to the convolution of arbitrarily dependent random variables, see \cite[Prop. 2.1]{cai:tang:2004} or \cite[Prop. 3.8(i)]{leipus:siaulys:konstantinides:2023}.

The Definition \ref{rem.KP.5.1}, permits arbitrarily dependent components of ${\bf X}$. If the following relation holds
\beam \label{eq.KP.5.9}
	\lim_{x\rightarrow\infty}\dfrac{P(X_{1}>x,\,\ldots,\,X_{n}>x)}{\overline{F}(x)}>0\,,
\eeam
then the components of $\bf{X}$ are pairwise asymptotically dependent. Furthermore, relation \eqref{eq.KP.5.2} is equivalent to
 \beao
 	\liminf_{\xto}\dfrac{\PP[X_1>b_1\,t_1\,x,\,\ldots,\,X_n>b_n\,t_n\,x]}{\PP[X_1>t_1\,x,\,\ldots,\,X_n>t_n\,x]} >0\,,
 \eeao
for any ${\bf t} \in (0,\,\infty]^n \setminus \{(\infty,\,\ldots,\,\infty)\}$ and any (or equivalently, for some) ${\bf b} > {\bf 1}$ (namely $b_{i}>1$ for any $i=1,\,\ldots,\,n$). We observe that if the components of the random vector ${\bf X}$ are independent and have weak equivalent marginal tails, then we get ${\bf F} \in \mathcal{D}_n$. Furthermore, we find out that because the definition of $\mathcal{D}_n$ requires that relation \eqref{eq.KP.5.2}  holds for any ${\bf t} \in (0,\,\infty]^n \setminus \{(\infty,\,\ldots,\,\infty)\}$, and having in mind the reductions of the dimension, when appears  $t_i=\infty$, for $i=1,\,\ldots,\,n$, it follows the implication: If  ${\bf F} \in \mathcal{D}_n$, then  ${\bf F}_{(n-k)} \in \mathcal{D}_{n-k}$, for any $k=1,\,\ldots,\,n-1$, where by ${\bf F}_{(n-k)}$ is denoted the distribution of the vector, that contains $n-k$ components from the initial $n$ ones (namely, reducing each time the $k$ components with infinite value for the components of ${\bf t}$). 
\ere

\bco \label{cor.KP.6.0}
Let  ${\bf X}$ be random vector, with distribution ${\bf F} \in \mathcal{D}_n$. Then it holds
\beam \label{eq.KP.6.a}
\PP[X_1>x,\,\ldots,\,X_n>x] \asymp \PP[X_1>t_1\,x,\,\ldots,\,X_n>t_n\,x]\,,
\eeam
as $\xto$, for any ${\bf t} \in (0,\,\infty]^n \setminus \{(\infty,\,\ldots,\,\infty)\}$. 
\eco

\pr~
The relation \eqref{eq.KP.6.a} is equivalent to
\beam \label{eq.KP.6.c}
0<\liminf_{\xto} \dfrac{\PP[X_1>t_1\,x,\,\ldots,\,X_n>t_n\,x]}{\PP[X_1>x,\,\ldots,\,X_n>x]} \leq \limsup_{\xto} \dfrac{\PP[X_1>t_1\,x,\,\ldots,\,X_n>t_n\,x]}{\PP[X_1>x,\,\ldots,\,X_n>x]} <\infty \,,
\eeam

We study separately three cases. In the first case, we consider ${\bf t} \in (0,\,1]^n$, and since from the property of class $\mathcal{D}_n$ we get directly the upper bound in \eqref{eq.KP.6.c},  for the lower bound we find
\beao
\liminf_{\xto} \dfrac{\PP[X_1>t_1\,x,\,\ldots,\,X_n>t_n\,x]}{\PP[X_1>x,\,\ldots,\,X_n>x]} \geq \liminf_{\xto} \dfrac{\PP[X_1>x,\,\ldots,\,X_n>x]}{\PP[X_1>x,\,\ldots,\,X_n>x]}=1>0\,,
\eeao 
for any ${\bf t} \in (0,\,1]^n$, hence we established \eqref{eq.KP.6.c}.

In the second case, we examine ${\bf t} \in (1,\,\infty]^n \setminus \{(\infty,\,\ldots,\,\infty)\}$. Now the lower bound of \eqref{eq.KP.6.c} is provided directly by class $\mathcal{D}_n$ properties. So, it remains  to check only if the upper bound is finite. Indeed, it holds
\beao
\limsup_{\xto} \dfrac{\PP[X_1>t_1\,x,\,\ldots,\,X_n>t_n\,x]}{\PP[X_1>x,\,\ldots,\,X_n>x]}\leq \limsup_{\xto} \dfrac{\PP[X_1>x,\,\ldots,\,X_n>x]}{\PP[X_1>x,\,\ldots,\,X_n>x]}=1<\infty\,.
\eeao
Here, we should notice that in case $t_i=\infty$, for $i=1,\,\ldots,\,n$, we find only dimension reduction. 

Finally, in the third case, whence ${\bf t} = (t_1,\,\ldots,\,t_n)$ has some components in the interval $t_i \in (0,\,1]$ and the others in the interval $t_j \in (1,\,\infty]$, for at least one $i\neq j \in \{1,\,\ldots,\,n\}$. Without loss of generality, for the sake of simplicity, we consider the case $t_1 \in (0,\,1]$ and $t_j \in (1,\,\infty]$ for $j=2,\,\ldots,\,n$, and all the other cases follow similarly. Thus, for the upper bound we obtain
\beao
&&\limsup_{\xto}\dfrac{\PP[X_1> t_1\,x,\,X_2>t_2\,x,\,\ldots,\,X_n>t_n\,x]}{\PP[X_1> x,\,X_2>x,\,\ldots,\,X_n>x]} \\[2mm]
&&\leq \limsup_{\xto}\dfrac{\PP[X_1> t_1\,x,\,X_2>x,\,\ldots,\,X_n>x]}{\PP[X_1> x,\,X_2>x,\,\ldots,\,X_n>x]} <\infty\,,
\eeao
where in the last step we used the  class $\mathcal{D}_n$ property of ${\bf X}$, since $t_1 \in (0,\,1]$. For the lower bound we obtain
\beao
&&\liminf_{\xto}\dfrac{\PP[X_1> t_1\,x,\,X_2>t_2\,x,\,\ldots,\,X_n>t_n\,x]}{\PP[X_1> x,\,X_2>x,\,\ldots,\,X_n>x]} \\[2mm]
&&\geq \liminf_{\xto}\dfrac{\PP[X_1> x,\,X_2>t_2\,x,\,\ldots,\,X_n>t_n\,x]}{\PP[X_1> x,\,X_2>x,\,\ldots,\,X_n>x]} > 0\,,
\eeao
where in the last step we used again the $\mathcal{D}_n$ property for ${\bf X}$,  and the fact that $t_j \in (1,\,\infty]$ for $j=2,\,\ldots,\,n$ and this completes the proof. 
~\halmos

From Corollary \ref{cor.KP.6.0}, follows that, if ${\bf X}$ has distribution ${\bf F} \in \mathcal{D}_n$, and satisfies the relation
 \beam \label{eq.KP.6.d}
\lim_{\xto} \dfrac{\PP[X_1>t_1\,x,\,\ldots,\,X_n>t_n\,x]}{\bF(x)} >0 \,,
\eeam
for some ${\bf t} \in (0,\,\infty]^n \setminus \{(\infty,\,\ldots,\,\infty)\}$, then has pairwise asymptotically dependent components. This permits even more flexibility in comparison with \eqref{eq.KP.5.9}.

\bre \label{rem.KP.6.1,5}
Let assume that $F \in \mathcal{D}$ and 
\beao
\bF_i(x) \asymp \bF(x)\,,
\eeao 
as $\xto$, where $F_i$ denotes a marginal distributions of vector ${\bf X}$, for $i=1,\,\ldots,\,n$, and let assume the asymptotic dependence among components, namely relation \eqref{eq.KP.6.d} holds, for all ${\bf t} \in (0,\,\infty]^n \setminus \{(\infty,\,\ldots,\,\infty)\}$. Then these conditions are sufficient to obtain  ${\bf F} \in \mathcal{D}_n$, that means to get \eqref{eq.KP.5.2}. Indeed, for any ${\bf t} \in (0,\,\infty]^n \setminus \{(\infty,\,\ldots,\,\infty)\}$ and for any ${\bf b} \in (0,\,1)^n$ we find
\beao
\limsup_{\xto} \dfrac{\PP[X_1>b_1\,t_1\,x,\,\ldots,\,X_n>b_n\,t_n\,x]}{\PP[X_1>t_1\,x,\,\ldots,\,X_n>t_n\,x]}\leq \limsup_{\xto} \dfrac{\delta\,\PP[X_1>t_1\,x]}{\alpha\,\bF(x)}\leq \dfrac {\delta}{\alpha}\,C_1< \infty\,,
\eeao
where the real number $\alpha\in (0,\,1]$ follows from the asymptotic dependence among the components of ${\bf X}$, the real number $0<\delta < \infty$ follows from the fact that these components belong to $\mathcal{D}_n$ and $C_1>0$ appears because of 
\beao
\bF_i(x) \asymp \bF(x)\,,
\eeao 
as $\xto$, for any $i=1,\,\ldots,\,n$, and by properties of class $\mathcal{D}$. The same happens  when the components are independent, and the proof is straight-forward. Hence, in these two extreme cases of dependence, namely the complete independence and the asymptotic dependence, relation \eqref{eq.KP.5.2} is unnecessary for the definition of class $\mathcal{D}_n$, however there exist several intermediate cases of dependence, that do not imply \eqref{eq.KP.5.2} directly.
\ere

Next we consider closure properties for $\mathcal{D}_n$ with respect to the product distribution. Here, for the product distribution of $Y$ and ${\bf X}$, we keep in mind the distribution of their scalar multiplication, namely 
\beao
{\bf \bH}({\bf t},\,x):=\PP[Y\,{\bf X}> {\bf t}\,x]=\PP[Y\,X_1>t_1\,x,\,\ldots,\,Y\,X_n>t_n\,x]\,,
\eeao 
for any  ${\bf t} \in (0,\,\infty]^n \setminus \{(\infty,\,\ldots,\,\infty)\}$. Therefore, assuming ${\bf F} \in \mathcal{D}_n$, in the case of independence between  $Y$ and ${\bf X}$, we wonder whether ${\bf H} \in \mathcal{D}_n$.

An application of the scalar products in insurance or financial framework is the following. Let us consider n-lines of business, the net-losses, described by the random variables $X_{1},\,\ldots,\,X_{n}$, and the random variable $Y$ denotes the common discount factor, over a given time interval. Then, the product $Y\,{\bf X}$ represents the discount net claims of the $n$ portfolios.

\begin{assumption} \label{ass.KP.5.1}
There exists function $a \;:\;[0\, ,\infty) \longrightarrow [0,\,\infty)$ such that hold
\begin{enumerate}

\item
$a(x) \longrightarrow \infty\,$,

\item
$a(x)=o(x)\,$, 

\item
$\bG[a(x)]=o\left[{\bf \bH}({\bf t},\,x)\right]\,$,
\end{enumerate}
as $\xto$, for any  ${\bf t} \in (0,\,\infty]^n \setminus \{(\infty,\,\ldots,\,\infty)\}$, where $G(x)=\PP[Y\leq x]$.
\end{assumption}

It is easy to see that Assumption \ref{ass.KP.5.1} is more strict than the Assumption \ref{ass.KP.A}. This can be verified by Remark \ref{rem.KP.1} and the inequalities 
\beao
{\bf \bH}({\bf 1},\,x) \leq \overline{H}_i(x)\,,
\eeao
for every $i=1,\,\ldots,\,n$. Because of the fact that the marginals of the vectors we use, follow heavy-tailed distributions, and therefore their right endpoint is infinite,  Assumption \ref{ass.KP.5.1}(3) is satisfied when distribution $G$ has bounded from above support.

\bth \label{th.KP.5.1}
Let ${\bf X}$ be random vector and  $Y$ be random variable, with distributions ${\bf F}$ and $G$ respectively, and let assume that $G(0-)=0\,,\; G(0)<1$. If  $Y$ and ${\bf X}$ are independent, Assumption  \ref{ass.KP.5.1} is valid and holds ${\bf F} \in \mathcal{D}_n$, then it holds ${\bf H} \in \mathcal{D}_n$. 
\ethe

\pr~
Likely to uni-variate case, we show only the case where $\bG(x)>0$, for any $x>0$, since in case $G$ has support with finite right endpoint, the argument follows similarly, putting instead of $a(x)$, the right endpoint of the support of $G$. 

Let ${\bf b} \in (0,\,1)^n$ and  ${\bf t} \in (0,\,\infty]^n \setminus \{(\infty,\,\ldots,\,\infty)\}$. From Assumption  \ref{ass.KP.5.1} we  obtain
\beam \label{eq.KP.5.4}
&&\limsup_{\xto}\dfrac {{\bf \bH}({\bf b}\,{\bf t},\,x)}{{\bf \bH}({\bf t},\,x)}=\limsup_{\xto}\dfrac{\PP[Y\,X_1>b_1\,t_1\,x,\,\ldots,\,Y\,X_n>b_n\,t_n\,x]}{\PP[Y\,X_1>t_1\,x,\,\ldots,\,Y\,X_n>t_n\,x]}\\[4mm] \notag
&&=\limsup_{\xto}\left(\int_0^{a(x)}+\int_{a(x)}^\infty \right)\dfrac{\PP\left[X_1>b_1\,t_1\,x/y,\,\ldots,\,X_n>b_n\,t_n\,x/y\right]}{\PP[Y\,X_1>t_1\,x,\,\ldots,\,Y\,X_n>t_n\,x]}\,G(dy)\\[4mm] \notag
&&=:\limsup_{\xto}\dfrac{I_1({\bf b}\,{\bf t},\,x)+I_2({\bf b}\,{\bf t},\,x)}{\PP[Y\,X_1>t_1\,x,\,\ldots,\,Y\,X_n>t_n\,x]}\,.
\eeam
Now, we calculate
\beao
I_2({\bf b}\,{\bf t},\,x) \leq \int_{a(x)}^\infty \,G(dy)=\bG[a(x)] = o\left[{\bf \bH}({\bf t},\,x)\right]\,,
\eeao
as $\xto$, due to Assumption  \ref{ass.KP.5.1}. Hence, taking into account relation \eqref{eq.KP.5.4} we get
\beao
\limsup_{\xto}\dfrac {{\bf \bH}({\bf b}\,{\bf t},\,x)}{{\bf \bH}({\bf t},\,x)}&\leq&\limsup_{\xto}\int_0^{a(x)}\dfrac{\PP\left[X_1>b_1\,t_1\,x/y,\,\ldots,\,X_n>b_n\,t_n\,x/y\right]}{\int_0^{a(x)}\PP\left[X_1>t_1\,x/y,\,\ldots,\,X_n>t_n\,x/y\right]\,G(dy)}\,G(dy)\\[2mm] 
&\leq&\limsup_{\xto} \sup_{0< y \leq a(x)} \dfrac{\PP\left[X_1>b_1\,t_1\,x/y,\,\ldots,\,X_n>b_n\,t_n\,x/y\right]}{\PP\left[X_1>t_1\,x/y,\,\ldots,\,X_n>t_n\,x/y\right]}\\[2mm] 
&=&\limsup_{\xto}\dfrac{\PP\left[X_1>b_1\,t_1\,x,\,\ldots,\,X_n>b_n\,t_n\,x\right]}{\PP\left[X_1>t_1\,x,\,\ldots,\,X_n>t_n\,x\right]}<\infty\,,
\eeao
where in the last step we used the condition ${\bf F} \in \mathcal{D}_n$. This way, we obtain that the multivariate tail, satisfies the desired property of $\mathcal{D}_n$. 

It remains now, to check if ${\bf H}$ has weak kernel. Indeed, it is enough to choose in the initial vector ${\bf F}$ as weak kernel $F$ some of the marginal distributions $F_i$, for $i=1,\,\ldots,\,n$. Then holds $\bF_i(x) \asymp \bF_j (x)$, as $\xto$, for any $i,\,j = 1,\,\ldots,\,n$, by definition of class $\mathcal{D}_n$, for ${\bf F}$ and further by Assumption \ref{ass.KP.5.1}, the inequality ${\bf \bH}({\bf 1},\,x) \leq \overline{H}_i(x)$ and the class $\mathcal{D}$ property, we obtain
\beao
\overline{H}_i(x)&=&\PP[Y\,X_i >x] =\left(\int_0^{a(x)} + \int_{a(x)}^\infty \right) \PP\left[ X_i >\dfrac xy \right]\,G(dy) \\[2mm]
&\leq& \int_0^{a(x)}\PP\left[ X_i >\dfrac xy \right]\,G(dy) +\bG[a(x)] \leq \PP\left[ X_i >\dfrac x{a(x)} \right] + o\left(\PP[Y\,X_i >x] \right)\,,
\eeao  
as $\xto$. Hence, from class $\mathcal{D}$ property and the previous relation, it follows
\beam \label{eq.KP.3.b}
\limsup_{\xto} \dfrac{\bH_i(x)}{\PP[X_i>x]} < \infty\,.
\eeam

From the other side, for any arbitrarily chosen $\vep \in (0,\,1)$, we have the convergence
\beam \label{eq.KP.3.c} \notag
\overline{H}_i(x)&\geq& \left(\int_{\vep}^1 + \int_1^{\infty} \right)\PP\left[X_i > \dfrac{x}{y}\right] \,\PP[Y \in dy] \\[2mm]
&\geq& d\,\PP[X_i>x] \,\PP[Y \in (\vep, \,1]] + \PP[X_i>x] \,\PP[Y > 1]\\[2mm] \notag
&\geq& (1\wedge d)\,\PP[X_i>x] \,\left(\PP[Y \in (\vep, \,1]] + \PP[Y > 1]\right) \longrightarrow (1\wedge d)\,\PP[X_i>x]\,,
\eeam 
as $\vep \downarrow 0$, where was used that $\PP\left[X_i > x/y \right]\geq \PP[X_i>x]$ in the second integral, that means when $y \in (1,\,\infty)$, while $d \in (0,\,1]$ comes as lower bound of this fraction for any $y \in (\vep, \,1)$, with $X_i$ from class $\mathcal{D}$. Thence by relation \eqref{eq.KP.3.c} we get
\beao
\liminf_{x \to \infty} \dfrac{\overline{H}_i(x)}{\PP[X_i>x]} \geq\liminf_{x \to \infty} \dfrac{(1\wedge d)\,\PP[X_i>x]}{\PP[X_i>x]} >0 \,,
\eeao
so we have
\beam \label{eq.KP.3.d}
\limsup_{x \to \infty} \dfrac{\PP[X_i>x] }{\overline{H}_i(x)} < \infty\,.
\eeam
Whence, from \eqref{eq.KP.3.b} and \eqref{eq.KP.3.d} follows the relation $\overline{H}_i(x) \asymp \PP[X_i >x]$, as $\xto$,  for any $i = 1,\,\ldots,\,n$ and therefore we get $\overline{H}_i(x) \asymp \overline{H}_j(x) $ as $\xto$, for any $i,\,j = 1,\,\ldots,\,n$, because $\bF_i(x) \asymp \bF_j(x)$ as $\xto$, for any $i,\,j = 1,\,\ldots,\,n$. Hence,  we get  $H_i \in \mathcal{D}$ by the closure property of class $\mathcal{D}$ with respect to weak equivalence (see \cite[Prop. 3.7(i)]{leipus:siaulys:konstantinides:2023}). Therefore, we find that ${\bf H}$ has weak kernel, and concretely, we can choose anyone from $H_1,\,\ldots,\,H_n$ as weak kernel. So we conclude that ${\bf H} \in \mathcal{D}_n$.
~\halmos

Let us introduce the following random vectors
\begin{enumerate}

\item
${\bf X}_1$ with distribution tail ${\bf \bF}_{1}({\bf t},\,x)=\PP[X_{11}> t_{1}\,x,\,\ldots,\,X_{1n} > t_{n}\,x]$, for some vector ${\bf t} \in (0,\,\infty]^n \setminus \{(\infty,\,\ldots,\,\infty)\}$.

\item
${\bf X}_2$ with distribution tail ${\bf \bF}_{2}({\bf t},\,x)=\PP[X_{21}> t_{1}\,x,\,\ldots,\,X_{2n} > t_{n}\,x]$, for some vector ${\bf t} \in (0,\,\infty]^n \setminus \{(\infty,\,\ldots,\,\infty)\}$.
\end{enumerate}
The sum of these two vectors has distribution tail $\overline{{\bf F}_{{\bf X}_1 +{\bf X}_2}}({\bf t},\,x):= \PP[{\bf X}_{1}+ {\bf X}_{2}> {\bf t}\,x]=\PP[X_{11}+X_{21}>t_{1}\,x,\,\ldots,\,X_{1n}+X_{2n} > t_{n}\,x]$, for some ${\bf t} \in (0,\,\infty]^n \setminus \{(\infty,\,\ldots,\,\infty)\}$ (when the random vectors ${\bf X}_1$ and ${\bf X}_2$ are independent then we denote $\overline{{\bf F}_{1} *{\bf F}_{2}}$  instead of ${\overline{{\bf F}}_{{\bf X}_1 +{\bf X}_2}}$). It is interesting to examine the closure property of the class $\mathcal{D}_{n}$ with respect to summation. For this reason we propose the following assumption. 

\begin{assumption}\label{ass.KP.5.3}
Any n-dimensional vector, consisting of some components of ${\bf X}_{1}$ and some components of ${\bf X}_{2}$, has distribution, which belongs to class $\mathcal{D}_{n}$.
\end{assumption}

The following example satisfies Assumption \ref{ass.KP.5.3}.

\bexam \label{exam.KP.5.3}
Let ${\bf X_{1}}$ and  ${\bf X_{2}}$ be independent random vectors, and hold ${\bf F}_1,\,{\bf F}_2 \in \mathcal{D}_n$, with $F_1,\,F_2$ their weak kernels respectively, satisfying $\bF_1(x)\asymp\bF_2(x)$, as $\xto$. For the sake of simplicity, let us consider the case $n=2$, for a vector whose one component is component of ${\bf X}_{1}$, and the other is component of ${\bf X}_{2}$, then the following upper limit is finite 
\beao
\limsup_{\xto}\dfrac{\PP[X_{11}>b_{1}\,t_{1}\,x,\; X_{22}>b_{2}\,t_{2}\,x]}{\PP[X_{11}>t_{1}\,x,\; X_{22}>t_{2}\,x]}
=	\limsup_{\xto}\dfrac{\PP[X_{11}>b_{1}\,t_{1}\,x]\,\PP[X_{22}>b_{2}\,t_{2}\,x]}{\PP[X_{11}>t_{1}\,x]\, \PP[X_{22}>t_{2}\,x]}\,,
\eeao
for every $b_{1},\,b_{2} \in (0,\,1)$,  because the distributions of $X_{11}$ and $X_{22}$ have distributions, which belong to class $\mathcal{D}$, since ${\bf F}_1,\,{\bf F}_2 \in \mathcal{D}_2$. Further, we can chose as weak kernel of the new random vector, any of the $F_1,\,F_2$ weak kernels, because of the assumption of weak equivalence between the two weak kernels. It is obvious that if one random vector have two components from the same primary vector (namely from ${\bf X}_{1}$ or from ${\bf X}_{2}$) then it belongs to class $\mathcal{D}_{2}$.
\eexam

Now we need a preliminary lemma, which represents partial generalization of \cite[Prop. 2.7]{shimura:watanabe:2005}, under arbitrary dependence but under the restriction in the frame of class $\mathcal{D}$. Let us remind, for any two non-negative random variables $X_1$ and $X_2$, are true the following elementary inequalities for $x\geq 0$
\beam \label{eq.KP.6.4}
\PP[{X}_{1}+ {X}_{2}>x] \leq \PP\left[X_{1}>\dfrac x2\right]+\PP\left[X_{2} >\dfrac x2\right]\,,\\[2mm] \label{eq.KP.6.5}
\PP[{X}_{1}+ {X}_{2}>x] \geq \dfrac 12 \,\left(\PP[X_{1}>x]+\PP[X_{2} >x]\right)\,.
\eeam

\ble \label{lem.KP.B*}
Let $X,\,X_i,\,Y,\,Y_i$ be non-negative random variables with the distributions $F$, $F_i$, $G$, $G_i \in \mathcal{D}$ respectively. If $\bF_i(x) \asymp \bF(x)$ and $\bG_i(x) \asymp \bG(x)$, as $\xto$, then it holds $\bF_{X_i+Y_i}(x) \asymp \bF_{X+Y}(x)$ as $\xto$, where $\bF_{X+Y}(x)=\PP[X+Y >x]$ and $\bF_{X_i+Y_i}(x)=\PP[X_i+Y_i >x]$. 
\ele

\pr~
Since $\bF_i(x) \asymp \bF(x)$ and $\bG_i(x) \asymp \bG(x)$, as $\xto$, we obtain
\beam \label{eq.KP.5.l1}
0 < \liminf_{\xto} \dfrac{\bF_i(x)}{\bF(x)} \leq \limsup_{\xto} \dfrac{\bF_i(x)}{\bF(x)} < \infty\,,\;\; 
0 < \liminf_{\xto} \dfrac{\bG_i(x)}{\bG(x)} \leq \limsup_{\xto} \dfrac{\bG_i(x)}{\bG(x)} < \infty\,,
\eeam
and taking into consideration that the random variables are non-negative, we apply the elementary inequalities \eqref{eq.KP.6.4} and \eqref{eq.KP.6.5} to find
\beao
\limsup_{\xto} \dfrac{\bF_{X_i+Y_i}(x)}{\bF_{X+Y}(x)} &\leq& 2\,\limsup_{\xto} \dfrac{\bF_{i}\left(\dfrac x2 \right) + \bG_{i}\left(\dfrac x2 \right)}{\bF(x)+\bG(x)}\leq 2\,\left(\limsup_{\xto} \dfrac{\bF_{i}\left(\dfrac x2 \right)}{\bF(x)}\bigvee  \dfrac{\bG_{i}\left(\dfrac x2 \right)}{\bG(x)}\right)  \\[2mm]
&<& \infty\,,
\eeao
where in the last step we used that $F_i,\,G_i \in \mathcal{D}$ and the relations \eqref{eq.KP.5.l1}. From the other side, we get similarly
\beao
\liminf_{\xto} \dfrac{\bF_{X_i+Y_i}(x)}{\bF_{X+Y}(x)} &\geq& \dfrac 12 \,\liminf_{\xto} \dfrac{\bF_{i}\left(x \right) + \bG_{i}\left(x \right)}{\bF\left(\dfrac x2 \right)+\bG\left(\dfrac x2 \right)} \\[2mm] \notag
&\geq& \dfrac 12 \,\left(\liminf_{\xto} \dfrac{\bF_{i}\left(x \right)}{\bF\left(\dfrac x2 \right)}\bigwedge \liminf_{\xto} \dfrac{\bG(x)}{\bG\left(\dfrac x2 \right)} \right)> 0\,.
\eeao
From the previous inequalities we conclude that $\bF_{X_i+Y_i}(x) \asymp \bF_{X+Y}(x)$, as $\xto$, for any $i=1,\,\ldots,\,n$.
~\halmos
 
\bth \label{th.KP.6.3}
Let ${\bf X}_1,\;{\bf X}_2$ be non-negative, arbitrarily dependent random vectors, with distributions ${\bf F}_{1},\;{\bf F}_{2} \in \mathcal{D}_n$ respectively, and assume Assumption \ref{ass.KP.5.3} is true. Then it holds 
${\bf F}_{{\bf X}_1 +{\bf X}_2} \in \mathcal{D}_n$. 
\ethe

\pr~
At first, we show the assertion for the case of $n=2$ and next we apply induction. 

Indeed for $n=2$ we obtain $\PP[{\bf X}_{1}+ {\bf X}_{2}> {\bf t}\,x]=\PP[X_{11}+X_{21}>t_{1}\,x,\,X_{12}+X_{22} >t_{2}\,x]$. So, by \eqref{eq.KP.6.4} we find
\beam  \label{eq.KP.6.6} \notag
&&\PP[{\bf X}_{1}+ {\bf X}_{2}> {\bf t}\,x] = \PP[X_{11}+X_{21}>t_{1}\,x\,,\;X_{12}+X_{22} >t_{2}\,x]\\[2mm] \notag
&&\leq \PP\left[X_{11}>\dfrac {t_1}2\,x\,,\;X_{12}+X_{22}>\dfrac {t_2}2\,x\right] +\PP\left[X_{21}>\dfrac {t_1}2\,x\,,\;X_{12}+X_{22}>\dfrac {t_2}2\,x\right]\\[2mm] 
&&\leq \PP\left[X_{11}>\dfrac {t_1}2\,x\,,\;X_{12}>\dfrac {t_2}2\,x\right]  +\PP\left[X_{11}>t_{1}2\,x\,,\;X_{22}> t_{2}2\,x\right]\\[2mm] \notag
&&\; +\PP\left[X_{21}>\dfrac {t_1}2\,x\,,\;X_{12}>\dfrac {t_2}2\,x\right]+\PP\left[X_{21}>\dfrac {t_1}2\,x\,,\;X_{22}>\dfrac {t_2}2\,x\right]\,.
\eeam
For the lower bound, we use successively \eqref{eq.KP.6.5} to find
\beam \label{eq.KP.6.7} \notag
&&\PP[{\bf X}_{1}+ {\bf X}_{2}> {\bf t}\,x] \\[2mm]  \notag
&& \geq \dfrac 12 \left( \PP[X_{11}>t_{1}\,x\,,\;X_{12}+X_{22} >t_{2}\,x]  +\PP\left[X_{21}>t_{1}\,x\,,\;X_{12}+X_{22}>t_{2}\,x\right]\right) \\[2mm]
&&\geq \dfrac 14 \Big( \PP[X_{11}>t_{1}\,x\,,\;X_{12}>t_{2}\,x] +\PP\left[X_{11}>t_{1}x\,,\;X_{22}>t_{2}\,x\right]  \\[2mm]  \notag
&&\qquad \qquad +\PP\left[X_{21}>t_{1}\,x\,,\;X_{12}>t_{2}\,x\right] +\PP\left[X_{21}>t_{1}\,x\,,\;X_{22}>t_{2}\,x\right] \Big)\,.
\eeam
In relations \eqref{eq.KP.6.6} and \eqref{eq.KP.6.7} there are joint probabilities, whose random variables belong to the same random vector, therefore they satisfy assumption $\mathcal{D}_2$, but there are some other joint probabilities whose random variables do not belong to the same random vector, as for example in the probability $\PP[X_{21}>t_{1}\,x\,,\;X_{12}> t_{2}\,x]$, but these probabilities also belong in the class $\mathcal{D}_2$ because of Assumption \ref{ass.KP.5.3}. Hence, from relations \eqref{eq.KP.6.6} and \eqref{eq.KP.6.7} we obtain the bound
\beao  
\limsup_{\xto}\dfrac{\PP\left[{\bf X}_{1}+ {\bf X}_{2}> {\bf t} \dfrac x2\right]}{\PP[{\bf X}_{1}+ {\bf X}_{2}> {\bf t}\,x]} \leq 4 \limsup_{\xto} \dfrac{P_1(x)+P_2(x)+P_3(x)+P_4(x)}{P'_1(x)+P'_2(x)+P'_3(x)+P'_4(x)}\leq 4 A^{2^2}(x)<\infty,
\eeao 
where
\beao
P_1(x)&=&\PP\left[X_{11}>t_{1}\,\dfrac x4\,,\;X_{12}>t_{2}\,\dfrac x4\right]\,,\quad P_2(x)=\PP\left[X_{11}>t_{1}\,\dfrac x4\,,\;X_{22}>t_{2}\,\dfrac x4\right]\,,\\[2mm]
P_3(x)&=&\PP\left[X_{21}>t_{1}\,\dfrac x4\,,\;X_{12}>t_{2}\,\dfrac x4\right]\,,\quad P_4(x)=\PP\left[X_{21}>t_{1}\,\dfrac x4\,,\;X_{22}>t_{2}\,\dfrac x4\right]\,,\\[2mm]
P'_1(x)&=&\PP[X_{11}>t_{1}\,x\,,\;X_{12}>t_{2}\,x]\,,\qquad P'_2(x)=\PP[X_{11}>t_{1}\,x\,,\;X_{22}>t_{2}\,x]\,,\\[2mm]
P'_3(x)&=&\PP[X_{21}>t_{1}\,x\,,\;X_{12}>t_{2}\,x]\,,\qquad P'_4(x)=\PP[X_{21}>t_{1}\,x\,,\;X_{22}>t_{2}\,x]\,,
\eeao
and $A^{2^2}(x)$ denotes the maximum of the upper limits of the fractions, as is shown in the expression 
\beao
&&A^{2^2}(x)\\[2mm] 
&&=\max\Bigg\{\limsup_{\xto}\dfrac{\PP\left[X_{11}>t_{1} \dfrac x4 ,\;X_{12}>t_{2} \dfrac x4 \right]}{\PP[ X_{11}>t_{1}\,x\,,\;X_{12}>t_{2}\,x]},\;\limsup_{\xto}\dfrac{\PP\left[X_{11}>t_{1} \dfrac x4,\,X_{22}>t_{2}\dfrac x4\right]}{\PP[ X_{11}>t_{1}\,x\,,\;X_{22}>t_{2}\,x]},\\[2mm] 
&&\limsup_{\xto}\dfrac{\PP\left[X_{21}>t_{1} \dfrac x4\,,\;X_{12}>t_{2} \dfrac x4\right]}{\PP\left[ X_{21}>t_{1}\,x\,,\;X_{12}>t_{2}\,x\right]},\;\limsup_{\xto}\dfrac{\PP\left[X_{21}>t_{1} \dfrac x4\,,\;X_{22}>t_{2} \dfrac x4 \right]}{\PP\left[ X_{21}>t_{1}\,x\,,\;X_{22}>t_{2}\,x\right]}\Bigg\}\,.
\eeao 
Because each term belongs in class $\mathcal{D}_2$, we find that $A^{2^2}(x)< \infty$, and furthermore by Lemma \ref{lem.KP.B*} we obtain that $
\bF_{X_1+X_2}(x) \asymp \bF_{X_{1i}+X_{2i}}(x)$, as $\xto$, where $F_{X_{1}}$, $F_{X_{2}}$ are the weak kernels of ${\bf F}_{1},\;{\bf F}_{2}$ respectively, therefore we established the relation ${\bf F}_{{\bf X}_1 +{\bf X}_2} \in \mathcal{D}_2$.

Next, we employ an induction argument. From the closure property with respect to sum for $n=2$, we assume that for some $n-1\geq 2$ the following inequalities hold
\beam \label{eq.KP.6.14}
\PP[{\bf X}_{1}' + {\bf X}_{2}'>{\bf t}\,x]\leq \sum_{z_i\in \bbz}\PP(z_i)\,, \qquad
\PP[{\bf X}_{1}' + {\bf X}_{2}'>{\bf t}\,x]\geq \dfrac 1{2^{n-1}}\sum_{z^*_i\in \bbz^*}\PP(z^*_i)\,,
\eeam 
where the ${\bf X}_{1}'$ and ${\bf X}_{2}'$ represent  $(n-1)$-dimensional vectors. The set $\bbz$ contains all the events of the type $\left\{X_{i_1,1}>t_{1}\,x/2,\,\ldots,\, X_{i_{n-1},n-1}>t_{n-1}\,x/2 \right\}$, where $i_1,\,\ldots,\,i_{n-1} \in \{1,\,2\}$. The set $\bbz^*$ contains all the events of the form $
\left\{X_{i_1,1}>t_{1}\,x,\,\ldots,\, X_{i_{n-1},n-1}>t_{n-1}\,x\right\}$, where $i_1,\,\ldots,\,i_{n-1} \in \{1,\,2\}$. We can observe that the  $\bbz$,  $\bbz^*$ have cardinal number $2^{n-1}$.

Now, we show the first relation in \eqref{eq.KP.6.14} with $n$, for ${\bf X}_1 \in  \mathcal{D}_n$ and ${\bf X}_2 \in  \mathcal{D}_n$.
\beam \label{eq.KP.6.16} \notag
\PP[{\bf X}_{1}+ {\bf X}_{2}> {\bf t}\,x]
&\leq& \PP\bigg[X_{11} + X_{21}>t_{1}\,x,\,\ldots,\, X_{1(n-1)} + X_{2(n-1)}>t_{n-1}\,x,\; X_{1n}>\dfrac{t_{n}}2\,x \bigg]\\[2mm]  \notag
&&+  \PP\bigg[X_{11} + X_{21}>t_{1}\,x,\,\ldots,\, X_{1(n-1)} + X_{2(n-1)}>t_{n-1}\,x,\; X_{2n}>\dfrac{t_{n}}2\,x \bigg]\\[2mm] 
&\leq& \sum_{z_i \in \bbz}\PP\left[z_i,\;X_{1n}>\dfrac{t_{n}}2\,x \right]+\sum_{z_i \in \bbz}\PP\left[z_i,\;X_{2n}>\dfrac{t_{n}}2\,x \right]\,,
\eeam
where the last step follows from the induction hypothesis in \eqref{eq.KP.6.14}. Therefore, from \eqref{eq.KP.6.16} we find that the probability $\PP[{\bf X}_{1}+ {\bf X}_{2}> {\bf t}\,x]$ is bounded from above by the sum of finite terms, whose cardinal number is $2^n$, and represent $n$-dimensional vectors.

Next, we consider the second relation in \eqref{eq.KP.6.14}. Through relation \eqref{eq.KP.6.5} and the inductive assumption in \eqref{eq.KP.6.14}, we find
\beam \label{eq.KP.6.17} \notag
\PP[{\bf X}_{1}+ {\bf X}_{2}> {\bf t}\,x]&\geq& \dfrac 12 \,\Big( \PP[X_{11} + X_{21}>t_{1}\,x,\,\ldots,\, X_{1(n-1)} + X_{2(n-1)} >t_{n-1}x,\, X_{1n}>t_{n}x ]\\[2mm]   \notag
&&+\PP[X_{11} + X_{21}>t_{1}\,x,\, \ldots,\, X_{1(n-1)} + X_{2(n-1)}>t_{n-1}\,x,\;X_{2n}>t_{n}\,x ] \Big)\\[2mm]  
&\geq& \dfrac 1{2^n} \,\left(\sum_{z^*_i \in Z^*} \PP\left[z^*_i,\;X_{1n}>t_{n}\,x\right] + \sum_{z^*_i \in Z^*} \PP\left[z^*_i,\; X_{2n}>t_{n}\,x\right] \right)\,.
\eeam
Thence, by relations \eqref{eq.KP.6.16} and \eqref{eq.KP.6.17}, we establish both inequalities \eqref{eq.KP.6.14} for any $n\geq 3$. 

Next, we see that all $n$-dimensional vectors from the right hand sides in relations \eqref{eq.KP.6.16} and \eqref{eq.KP.6.17} belong to class $\mathcal{D}_{n}$, and from these formulas we obtain
\beam \label{eq.KP.6.18} 
\limsup_{\xto}\dfrac{\PP\left[{\bf X}_{1}+ {\bf X}_{2} >{\bf t}\,\dfrac x2\right]}{\PP[{\bf X}_{1}+ {\bf X}_{2}> {\bf t}\,x]}\leq 2^n\,A^{2^n}(x) < \infty\,,
\eeam
where $A^{2^n}(x)$ is defined as maximum of the corresponding $2^n$ fractions, containing the distribution tails of $n$-dimensional random vectors, similarly to  $A^{2^2}(x)$. The $A^{2^n}(x)$ is finite, as far as all the $n$-dimensional random vectors belong to class $\mathcal{D}_{n}$ and again from Lemma \ref{lem.KP.B*} the components of the summands have weak equivalent tails, so anyone of the marinals of the new vector ${\bf X}_{1}+ {\bf X}_{2}$, can play the role of weak kernel, since the class $\mathcal{D}$ has the closure property with respect to sum. Hence  ${\bf F}_{{\bf X}_1 +{\bf X}_2} \in \mathcal{D}_n$ for any $n\in \bbn$.
~\halmos

\bre \label{rem.KP.6.1}
Relation \eqref{eq.KP.6.18} renders the closure property of class $\mathcal{D}_{n}$ under non-negative, arbitrarily dependent random vectors, each of which  has arbitrarily dependent components. This has important consequences, as in most risk models is allowed some dependence among the components of the random vectors but among the random vectors still remains independence. Furthermore, in most cases is applied the multivariate regular variation framework, which belongs to class $\mathcal{D}_{n}$, as it follows from Theorem \ref{th.KP.6.2} below, in case of asymptotic dependent components. For example, on these risk models we refer to \cite{konstantinides:li} and \cite{li:2016}. In some other works the dependence among random vectors but with their components mutually independent was studied, as for example in \cite{yang:li:2014} and \cite{yang:li:2015}.
\ere

Next we introduce the randomly stopped sum $S_{N}:=\sum_{i=1}^{N}X_{i}$, with $S_{0}=0$, where $\{X_{i},\; i\in \bbn\}$ represents a sequence of random variables, with distributions $\{F_{i},\; i\in \bbn\}$ respectively, and $N$ represents a discrete random variable, independent of $\{X_{i},\; i\in \bbn\}$, whose distribution has as support some subset of $\{0,1,...\}$. The probability mass function of $N$ can be given through $p_{n}:=P[N=n]$, with $p_{0}<1$ and 
\beao
\sum_{n=0}^{\infty}p_{n}=1\,.
\eeao 
We should mention that if $N$ has support with lower bound $0$, then put $\kappa = 1$, while if it has support with lower bound an integer greater or equal to unity, then we denote this bound by $\kappa$. The randomly stopped sum $S_{N}$ has distribution $F_{S_{N}}$, whose tail is given by 
\beao
\overline{F}_{S_{N}}(x):=P(S_{N}>x)=\sum_{n=1}^{\infty}p_{n}P(S_{n}>x)\,,
\eeao 
for any $x>0$, with 
\beao
S_{n}:=\sum_{i=1}^{n}X_{i}\,.
\eeao 
For example, in insurance framework the sequence of $\{X_{i},\; i\in \bbn\}$ represents the losses (or profits), and the random variable $N$ corresponds to the number of losses over the given time interval. As a result the quantity $S_{N}$ represents the total loss of the insurance business. For some papers related to $S_{N}$, in presence of heavy tails, we refer to \cite{karaseviciene:siaulys:2023}, \cite{leipus:siaulys:2012}, and \cite{watanabe:2008} among others.

Following the previous arguments, we examine some closure properties for class $\mathcal{D}_{n}$, in multivariate frame. Namely, we take into consideration the finite mixing property, as also the randomly stopped sum 
\beao
{\bf S}_{N}=\sum_{i=1}^N {\bf X}_{i}\,,
\eeao 
with ${\bf S_{0}}={\bf 0}$, and here $N$ is defined as above and the distribution tail of ${\bf S}_N$ takes the form 
\beao
{\bf \bF}_{{\bf S}_N}({\bf t},\,x)=\PP\left[{\bf S}_{N}>{\bf t} \,x\right]=\sum_{d=1}^{\infty} p_d\,\PP[{\bf S}_{d}>{\bf t}\,x]\,,
\eeao 
with $p_d=\PP[N=d]$ with $p_0<1$ and 
\beao
\sum_{j=0}^{\infty} p_j =1\,,
\eeao 
when $x>0$, where
\beam \label{eq.KP.6.21} 
{\bf \bF}_{{\bf S}_d}({\bf t},\,x)=\PP\left[{\bf S}_{d}>{\bf t} \,x\right]=\PP[X_{11}+\cdots +X_{d1}>t_{1}\,x,\,\ldots,\,X_{1n}+\cdots +X_{dn}>t_{n}\,x]\,,
\eeam
for every ${\bf t} \in (0,\,\infty]^n \setminus \{(\infty,\,\ldots,\,\infty)\}$.

\bre \label{rem.KP.A*}
Let ${\bf X}_{1}, \,{\bf X}_{2}$ be $n$-dimensional random vectors, with distributions ${\bf F}_{1},\,{\bf F}_{2}$ respectively. If we denote by ${\bf M}$ the finite mixture, and by $ \overline{\bf M} ({\bf t},\,x)=p\,{\bf \overline{F}}_{1}({\bf t},\,x)+(1-p)\,{\bf \overline{F}}_{2}({\bf t},\,x)$ its tail, for any $p \in (0,\,1)$, namely for any ${\bf t} \in (0,\,\infty]^n \setminus \{(\infty,\,\ldots,\,\infty)\}$, 
\beao
 \overline{\bf M} ({\bf t},\,x)=p\,\PP[X_{11}>t_1\,x,\,\ldots,\,X_{1n}>t_n\,x]+(1-p)\,\PP[X_{21}>t_1\,x,\,\ldots,\,X_{2n}>t_n\,x]\,,
\eeao
then we call as weak kernel of ${\bf M} $, some distribution $M$,  representing a finite mixture of the weak kernels $F_1,\,F_2$ of the ${\bf F}_{1},\,{\bf F}_{2}$ respectively, $M=p\,F_{1} +(1-p)\,F_{2}$, for any  $p \in (0,\,1)$.

We should remark that since the weak kernels of ${\bf F}_{1},\,{\bf F}_{2}$ satisfy the relations $\bF_1(x) \asymp \bF_{1i}(x)$ and $\bF_2(x) \asymp \bF_{2i}(x)$, as $\xto$, for any $i=1,\,\ldots,\,n$, where by $F_{ji}$ we denote the $i$-th marginal of $j$-th vector, with $j=1,\,2$, it follows that $\overline{M}(x) \asymp \overline{M}_i(x)$, as $\xto$, for any $i=1,\,\ldots,\,n$, with $M_i:=p\,F_{1i} + (1-p)\,F_{2i}$. Indeed, from one side, since the weak kernels satisfy the weak equivalence, then 
\beam \label{eq.KP.5.r4} \notag
\limsup_{\xto} \dfrac{\overline{M}_i(x)}{\overline{M}(x)} &=& \limsup_{\xto} \dfrac{p \bF_{1i}(x) + (1-p) \bF_{2i}(x)}{p\,\bF_{1} + (1-p)\,\bF_{2}(x)}\\[2mm]
&\leq& \limsup_{\xto} \left(\dfrac{\bF_{1i}(x)}{\bF_{1}(x)}\bigvee \dfrac{\bF_{2i}(x)}{\bF_{2}(x)}\right) < \infty\,,
\eeam
where in the last step we use that the weak kernels  $F_1,\,F_2$ are weak equivalent with respect to the corresponding marginals. 

From the other side, with similar arguments we find
\beam \label{eq.KP.5.r5} \notag
\liminf_{\xto} \dfrac{\overline{M}_i(x)}{\overline{M}(x)} &=& \liminf_{\xto} \dfrac{p\,\bF_{1i}(x) + (1-p)\,\bF_{2i}(x)}{p\,\bF_{1} + (1-p)\,\bF_{2}(x)}\\[2mm]
&\geq& \liminf_{\xto} \left( \dfrac{\bF_{1i}(x)}{\bF_{1}(x)}\bigwedge \liminf_{\xto} \dfrac{\bF_{2i}(x)}{\bF_{2}(x)} \right)> 0\,.
\eeam
So from relations \eqref{eq.KP.5.r4} and \eqref{eq.KP.5.r5} follows that if $\bF_1(x) \asymp \bF_{1i}(x)$ and $\bF_2(x) \asymp \bF_{2i}(x)$, as $\xto$, for any $i=1,\,\ldots,\,n$, then it holds $\overline{M}_i(x)\asymp\overline{M}(x)$, as $\xto$, for any $i=1,\,\ldots,\,n$.
\ere

The next result extends some known partial results, as it contains arbitrarily dependent $X_1,\,X_2,\,\ldots$ but our restriction is the non-negative random variables. For some relative results in the class $\mathcal{D}$ see \cite[sec. 4]{karaseviciene:siaulys:2023}. 

\ble \label{lem.KP.5.G*}
Let $X_1,\,X_2,\,\ldots$ be non-negative, arbitrarily dependent  random variables with distributions $F_1,\,F_2,\,\ldots$ respectively, where $F_i \in \mathcal{D}$, for any  $i \in \bbn$. If $N$ represents a discrete, non-degenerated to zero random variable and is independent of $X_1,\,X_2,\,\ldots$, then it holds $F_{S_N} \in \mathcal{D}$.
\ele

\pr~
Without loss of generality, we assume, for sake of simplicity, that $k=1$. In case $N=m=1$, the result is trivial. In case $N=m$, with integer $m\geq 2$, we obtain that
\beao
\bF_{S_N}(x)=p_1\,\bF_1(x) + \cdots + p_m\,\bF_{X_1 + \cdots +X_m}(x)\,,
\eeao 
for any $x>0$. Since $X_i$ follows distribution $F_i \in \mathcal{D}$ for any  $i \in \bbn$ and they are non-negative random variables, we find by \cite[Prop. 2.1]{cai:tang:2004} or \cite[Prop. 3.8]{leipus:siaulys:konstantinides:2023} that $F_{X_1+\cdots+X_d} \in \mathcal{D}$. Hence, for any $n\in \bbn$ and any $b\in (0,\,1)$, there exists finite constant $c_d<\infty$, such that for any large enough $x>0$, it holds $\PP[S_d>b\,x] \leq c_d\,\PP[S_d>x]$, therefore we obtain $\widehat{c}:=\sup\{c_1,\,c_2,\,\ldots \} < \infty$. So we find
\beao
\bF_{S_N}(b\,x) = \sum_{d=1}^{\infty} \PP[S_d> b\,x]\,\PP[N=d] \leq \widehat{c}\,\sum_{d=1}^{\infty} \PP[S_d > x]\,\PP[N=d]=\widehat{c}\,\bF_{S_N}(x)\,,
\eeao
for any $b\in (0,\,1)$ and any $x>0$. That means 
\beao
\limsup \dfrac{\bF_{S_N}(b\,x)}{\bF_{S_N}(x)}  < \infty\,,
\eeao 
which implies that $F_{S_N} \in \mathcal{D}$.
~\halmos 

Now we are ready to give the closure property of class $\mathcal{D}_{n}$ with respect to finite mixture and the multivariate randomly stopped sums. We note that as in case of sum of two random vectors, when we employed Assumption \ref{ass.KP.5.3}, here we repeat the same for the sum of random number of summands, but now pairwisely. Namely, the conditions of Assumption \ref{ass.KP.5.3}, are satisfied  by ${\bf X}_{1}$ with ${\bf X}_{2}$, by ${\bf X}_{1}+ {\bf X}_{2}$ with ${\bf X}_{3}$, by ${\bf X}_{1}+ {\bf X}_{2}+{\bf X}_{3}$ with ${\bf X}_{4}$, etc. 
 
\bco \label{cor.KP.6.1}
\begin{enumerate}

\item
Let ${\bf X}_{1}, \,{\bf X}_{2}$ be arbitrarily dependent random vectors with distributions ${\bf F}_1,\,{\bf F}_2 \in \mathcal{D}_{n}$ respectively. Then it holds ${\bf M}\in \mathcal{D}_{n}$. 

\item
Let ${\bf X}_{1}, \,{\bf X}_{2},\,\ldots$ be non-negative, arbitrarily dependent random vectors with distributions ${\bf F}_{ 1},\,{\bf F}_{ 2},\,\ldots\, \in \mathcal{D}_{n}$ respectively, and Assumption \ref{ass.KP.5.3} be true pairwisely. If $N$ is a discrete random variable, independent from random vectors ${\bf X}_{1}, \,{\bf X}_{2},\,\ldots$, then it holds ${\bf F}_{{\bf S}_N} \in \mathcal{D}_{n}$.
\end{enumerate}
\eco

\pr~
\begin{enumerate}

\item
For any $p \in (0,\,1)$, and any ${\bf t} \in (0,\,\infty]^n \setminus \{(\infty,\,\ldots,\,\infty)\}$ and any ${\bf b}\,\in (0,\,1)^n$ holds
\beao
&&\limsup_{\xto} \dfrac{p\,{\bf \bF}_1({\bf b}\,{\bf t},\,x) + (1-p)\,{\bf \bF}_2({\bf b}\,{\bf t},\,x)}{p\,{\bf \bF}_1({\bf t},\,x) + (1-p)\,{\bf \bF}_2({\bf t},\,x)}\\[2mm]
&&= \limsup_{\xto}\dfrac{p\,\PP_1({\bf b}\,{\bf t})+ (1-p)\,\PP_2({\bf b}\,{\bf t})}{p\,\PP[X_{11}>t_{1}\,x,\,\ldots,\,X_{1n}>t_{n}\,x ]+ (1-p)\,\PP[X_{21}>t_{1}\,x,\,\ldots,\,X_{2n}>t_{n}\,x ]}\\[2mm]
&&\leq\max \left\{ \limsup_{\xto} \dfrac{\PP[X_{11}>b_{1}\,t_{1}\,x,\,\ldots,\,X_{1n}>b_{n}\,t_{n}\,x ]}{\PP[X_{11}>t_{1}\,x,\,\ldots,\,X_{1n}>t_{n}\,x ]}\,,\;\right.\\[2mm]
&&\qquad \qquad \left.\limsup_{\xto} \dfrac{\PP[X_{21}>b_{1}\,t_{1}\,x,\,\ldots,\,X_{2n}>b_{n}\,t_{n}\,x ]}{\PP[X_{21}>t_{1}\,x,\,\ldots,\,X_{2n}>t_{n}\,x ]}\right\} < \infty \,,
\eeao
where 
\beao
\PP_1({\bf b}\,{\bf t})=\PP[X_{11}>b_{1}\,t_{1}\,x,\,\ldots,\,X_{1n}>b_{n}\,t_{n}\,x ]\,,\\[2mm]
\PP_2({\bf b}\,{\bf t})=\PP[X_{21}>b_{1}\,t_{1}\,x,\,\ldots,\,X_{2n}>b_{n}\,t_{n}\,x ]\,,
\eeao 
and the last step is due to the assumptions ${\bf F}_1,\,{\bf F}_2 \in \mathcal{D}_{n}$. Further, by the relations $F_{1i} \in \mathcal{D}$, $F_{2i} \in \mathcal{D}$, and the closure properties with respect to finite mixture (see \cite[Prop. 3.7.(iv)]{leipus:siaulys:konstantinides:2023}) we obtain $M_i \in \mathcal{D}$ for any $i=1,\,\ldots,\,n$. 

Next, by Remark \ref{rem.KP.A*} we know that there exists weak kernel $M$ such that $\overline{M}(x) \asymp \overline{M}_i(x)$, as $\xto$ for any  $i=1,\,\ldots,\,n$, therefore we get $M \in \mathcal{D}$. So we conclude $
p\,{\bf F}_1+ (1-p)\,{\bf F}_2  \in \mathcal{D}_{n}$.

\item
When $N=m=1$, the result follows directly. We consider the case $N=m$ for some integer $m\geq 2$. Initially for such $m \in \bbn$, we obtain
\beam \label{eq.KP.6.22} 
{\bf \bF}_{{\bf S}_N}({\bf t},\,x)=\PP\left[{\bf S}_{N}>{\bf t} \,x\right]=p_1\,\PP[{\bf X}_{1}>{\bf t}\,x]+\cdots+p_m\,\PP[{\bf X}_{1}+\cdots +{\bf X}_{m}>{\bf t}\,x] \,,
\eeam
for any $x>0$,  ${\bf t}\in (0,\,\infty]^n \setminus \{(\infty,\,\ldots,\,\infty)\}$. By Theorem \ref{th.KP.6.3} we find that each term in the right hand side of \eqref{eq.KP.6.22} belongs to class $\mathcal{D}_{n}$. Hence, for any $x>0$, any ${\bf t}\in (0,\,\infty]^n \setminus \{(\infty,\,\ldots,\,\infty)\}$ and any ${\bf b} \in (0,\,1)^n$ we obtain
\beam \label{eq.KP.5.aaa}
\PP[{\bf S}_d > {\bf b}\,{\bf t}\,x] \leq c_d\,\PP[{\bf S}_d > {\bf t}\,x]\,,
\eeam
for some finite constant $c_d\geq 1$. Thus, for $\widehat{c}:=\sup\{c_1,\,c_2,\,\ldots\} < \infty$ 
\beam  \label{eq.KP.5.bbb} \notag
{\bf F}_{{\bf S}_N}({\bf b}\,{\bf t},\,x)&=& \sum_{n=1}^\infty \PP[{\bf S}_d > {\bf b}\,{\bf t}\,x]\,\PP[N=d] \leq \sum_{n=1}^\infty c_d\,\PP[{\bf S}_d > {\bf t}\,x]\,\PP[N=d]\\[2mm]
&\leq& \widehat{c}\,\sum_{n=1}^\infty \PP[{\bf S}_d > {\bf t}\,x]\,\PP[N=d]=\widehat{c}\,{\bf F}_{{\bf S}_N}({\bf t},\,x)\,,
\eeam
with ${\bf b} \in (0,\,1)^n$, for any $x>0$ and any ${\bf t}\in (0,\,\infty]^n \setminus \{(\infty,\,\ldots,\,\infty)\}$, where in second step we used relation \eqref{eq.KP.5.aaa}. Therefore from relation \eqref{eq.KP.5.bbb} and the arbitrary choice of ${\bf t}$, $x>0$ and ${\bf b}$ we find
\beao
\limsup_{\xto} \dfrac{{\bf F}_{{\bf S}_N}({\bf b}\,{\bf t},\,x)}{{\bf F}_{{\bf S}_N}({\bf t},\,x)} < \infty\,,
\eeao
which represents the multivariate property we need.

Similarly we examine if the multivariate randomly stopped sum has weak kernel. If the random variable $N$ takes the value $d \in \{1,\,2,\,\ldots \}$, where $p_d>0$, then by relation \eqref{eq.KP.6.21} we find that the $X_{11},\,\ldots,\,X_{d1}$  follow the distributions $F_{11},\,\ldots,\,F_{d1}$, that belong to $d$-different random vectors ${\bf X}_1,\,\ldots,\, {\bf X}_d$ respectively, with corresponding distributions ${\bf F}_1,\,\ldots,\, {\bf F}_d \in \mathcal{D}_n$, therefore all these multivariate distributions have weak kernels $F_{1},\,\ldots,\,F_{d}$ respectively. Thus we see that $\bF_{i1}(x) \asymp \bF_{i}(x)$, as $\xto$, for any $i=1,\,\ldots,\,d$. So, if the distribution of the sum 
\beao
S_{d1}:=\sum_{i=1}^d X_{i1}\,,
\eeao 
is denoted by $F_{S_{d1}}$ and the distribution of the  sum 
\beao
S_d=\sum_{i=1}^d X_i\,,
\eeao 
is denoted by $F_{S_d}$, (where $X_i$ are the corresponding non-negative random variables, whose distributions has weak kernels $F_i$), then applying Lemma \ref{lem.KP.B*} we obtain the weak equivalence 
\beao
\bF_{S_{d1}}(x) \asymp \bF_{S_d}(x)\,,
\eeao 
as $\xto$. Similarly for the rest components of the vector we get 
\beam \label{eq.KP.5.1t}
\bF_{S_{di}}(x) \asymp \bF_{S_d}(x)\,,
\eeam 
as $\xto$, for any  $i=1,\,\ldots,\,n$.

From \cite[Prop. 2.1]{cai:tang:2004}, we obtain $F_{S_{di}} \in \mathcal{D}$ and from relation \eqref{eq.KP.5.1t} we see that $F_{S_{d}} \in \mathcal{D}$, see \cite[Prop. 3.7]{leipus:siaulys:konstantinides:2023}. Hence,
\beam \label{eq.KP.5.2t}
\bF_{S_{Ni}}(x) \asymp \bF_{S_N}(x)\,,
\eeam 
as $\xto$, for any  $i=1,\,\ldots,\,n$, where $S_{Ni}:=\sum_{j=1}^N X_{ji}$, that follows distribution $F_{S_{Ni}}$ for any   $i=1,\,\ldots,\,n$ and 
$S_{N}:=\sum_{j=1}^N X_{j}$, that follows distribution $F_{S_{N}}$, namely it is the randomly stopped sum of each component. So, since applying Lemma \ref{lem.KP.5.G*} we obtain $F_{S_{Ni}}\in \mathcal{D}$, from relation \eqref{eq.KP.5.2t} we find $F_{S_{N}}\in \mathcal{D}$. From this and relation \eqref{eq.KP.5.2t} follows that  $F_{S_{N}}$ is the weak kernel of this randomly stopped sum of random vectors. In combination with the behavior of the tail, that was found finite, we reach the desired result.
~\halmos
\end{enumerate}

\bre \label{rem.KP.5.6}
A specific example of ${\bf S}_{N}$ is the compound Poisson, when the $N$ represents a random variable with distribution Poisson and the random vectors ${\bf X}_{1},\,{\bf X}_{2},\,\ldots$ are independent and identically distributed, with $N$ independent of ${\bf X_{i}}$ for any $i\in \bbn$. It is easy to see that in this case, when the vectors ${\bf X}_{1},\,{\bf X}_{2},\,\ldots$ follow distribution ${\bf F} \in \mathcal{D}_n$, then the distribution of ${\bf S}_{N}$ belongs also in $\mathcal{D}_n$, as follows by Corollary \ref{cor.KP.6.1}, since Assumption \ref{ass.KP.5.3} is valid for independent and identically distributed vectors, because we have a special case in Example \ref{exam.KP.5.3}. The compound Poisson  plays crucial role in L\'{e}vy processes, see for example \cite{das:fasenhartmann:2023} and its references.
\ere

\section{Positively decreasing random vectors} \label{sec.KP.7}

Analogously to class $\mathcal{D}_n$, we define the corresponding class $\mathcal{P_D}_{n}$, as multivariate analogue of the positive decrease $\mathcal{P_D}$.

\bde \label{def.KP.7.1}
Let ${\bf X}$ be a random vector, with distribution ${\bf F}$ and marginal distributions $F_1,\,\ldots,\,F_n$. Let assume that there exists a distribution $F$, named strong kernel, such that $F \in \mathcal{P_D}$ and $\bF_i(x) \sim c_i\,\bF(x)$, as $\xto$, with $c_i\in (0,\,\infty)$, for any $i=1,\,\ldots,\,n$. If holds the inequality
\beam \label{eq.KP.7.1}
\limsup_{\xto} \dfrac{{\bf \bF}({\bf v}\,{\bf t},\,x)}{{\bf \bF}({\bf t},\,x)}=\limsup_{\xto} \dfrac{\PP[X_1 > v_1\,t_1\,x,\,\ldots,\,X_n > v_n\,t_n\,x]}{\PP[X_1 > t_1\,x,\,\ldots,\,X_n > t_n\,x]} < 1\,,
\eeam
for any (or equivalently, for some) ${\bf v}=(v_1,\,\ldots,\,v_n) > {\bf 1}=(1,\,\ldots,\,1)$ and for all ${\bf t} \in (0,\,\infty]^n \setminus \{(\infty,\,\ldots,\,\infty)\}$, then we say that the random vector has multivariate positively decreasing distribution, symbolically ${\bf F} \in \mathcal{P_D}_{n}$.
\ede
While in class $\mathcal{D}_n$, there was a weak kernel, namely a distribution $F\in \mathcal{D}$, such that $\bF_i(x) \asymp \bF(x)$, as $\xto$, for any $i=1,\,\ldots,\,n$, here we restrict ourselves to a strong kernel, namely to a distribution $F$, such that $\bF_i(x) \sim c_i\,\bF(x)$, as $\xto$, for $c_i\in (0,\,\infty)$,  for any $i=1,\,\ldots,\,n$. The reason for this, is that class $\mathcal{P_D}$ is closed with respect to strong equivalence, but not with respect to weak equivalence. Hence, $F_i \in \mathcal{P_D}$, for all marginals of vectors in $\mathcal{P_D}_n$.  

In next result we present a closure property of class $\mathcal{P_D}_n$ with respect to scalar product under the assumption of identical components of ${\bf X}$.  

\bth \label{th.KP.7.1}
Let ${\bf X}$ be random vector and $Y$ be random variable, with distributions ${\bf F}$ and $G$ respectively, and  assume $G(0-)=0\,,\; G(0)<1$. Additionally we assume that ${\bf X}$ have identically distributed components with common distribution $F$.
If ${\bf F} \in \mathcal{P_D}_n$, and ${\bf X}$ and $Y$ are independent, satisfying Assumption \ref{ass.KP.5.1}, then it holds ${\bf H} \in \mathcal{P_D}_n$. 
\ethe

\pr~
As in Theorem \ref{th.KP.5.1}, here we show only the case $\bG(x)>0$, for any $x>0$, since in opposite case, we follows the same proof with only change of $a(x)$ to the upper endpoint of the support of $G$.

Let ${\bf v}=(v_1,\,\ldots,\,v_n) >{\bf 1}$ and ${\bf t} \in (0,\,\infty]^n \setminus \{(\infty,\,\ldots,\,\infty)\}$. From Assumption \ref{ass.KP.5.1} we find
\beao
&&\limsup_{\xto}\dfrac {{\bf \bH}({\bf v}\,{\bf t},\,x)}{{\bf \bH}({\bf t},\,x)}=\limsup_{\xto} \dfrac{\PP\left[Y\,X_1 > v_1\,t_1\,x,\,\ldots,\,Y\,X_n > v_n\,t_n\,x\right]}{\PP[Y\,X_1 > t_1\,x,\,\ldots,\,Y\,X_n > t_n\,x]}\\[2mm]
&&= \limsup_{\xto} \left(\int_0^{a(x)} + \int_{a(x)}^\infty \right)\dfrac{ \PP\left[X_1>v_1\,t_1\,\dfrac xy,\,\ldots,\,X_n>v_n\,t_n\,\dfrac xy\right]}{ \PP\left[Y\,X_1>t_1\,x,\,\ldots,\,Y\,X_n> t_n\,x \right]}\,G(dy)\\[2mm] 
&&=:\limsup_{\xto}\dfrac{I_1(x)+I_2(x)}{ \PP\left[Y\,X_1>t_1\,x,\,\ldots,\,Y\,X_n> t_n\,x \right]}\,.
\eeao
Hence, 
\beao
I_2(x) \leq \bG[a(x)] =o\left[{\bf \bH}({\bf t},\,x)\right]\,,
\eeao 
as $\xto$, and therefore
\beao
\limsup_{\xto}\dfrac {{\bf \bH}({\bf v}\,{\bf t},\,x)}{{\bf \bH}({\bf t},\,x)}&\leq&\limsup_{\xto}\int_0^{a(x)}\dfrac{ \PP\left[X_1>v_1\,t_1\,\dfrac xy,\,\ldots,\,X_n>v_n\,t_n\,\dfrac xy\right]}{\int_0^{a(x)} \PP\left[X_1>t_1\,\dfrac xy,\,\ldots,\,X_n>t_n\,\dfrac xy\right]\,G(dy)}\,G(dy)\\[2mm] 
&\leq&\limsup_{\xto} \sup_{0< y \leq a(x)} \dfrac{ \PP\left[X_1>v_1\,t_1\,\dfrac xy,\,\ldots,\,X_n>v_n\,t_n\,\dfrac xy\right]}{\PP\left[X_1>t_1\,\dfrac xy,\,\ldots,\,X_n>t_n\,\dfrac xy\right]}\\[2mm] 
&=&\limsup_{\xto} \dfrac{ \PP\left[X_1>v_1\,t_1\,x,\,\ldots,\,X_n>v_n\,t_n\,x\right]}{\PP\left[X_1>t_1\,x,\,\ldots,\,X_n>t_n\,x\right]} < 1\,,
\eeao
where in the last step we used that ${\bf F} \in \mathcal{P_D}_n$.
Since ${\bf X}$ has identically distributed components, then also the products $Y\,X_i$ for $i=1,\,\ldots,\,n$ are identically distributed. So, by Assumption \ref{ass.KP.5.1} and the fact that ${\bf \bH}({\bf 1},\,x) \leq \overline{H}_i(x)$, 
we obtain Assumption \ref{ass.KP.A} (see Remark \ref{rem.KP.1}), and since the independence is included in Assumption \ref{ass.KP.B}, applying Theorem \ref{th.KP.3} we find that the distribution of $Y\,X_i$ is $H_i \in \mathcal{P_D}$. So any $H_i$ can play the role of the strong kernel. 
Thus we get ${\bf H} \in \mathcal{P_D}_n$.
~\halmos

The next result is related with the intersection of these two multivariate classes. To avoid misunderstanding, let us remind that a random vector ${\bf X}$, with distribution ${\bf F}$, belongs to the class $\mathcal{D}_n \cap \mathcal{P_D}_n $ if there exists strong kernel $F \in \mathcal{D} \cap \mathcal{P_D}$ and it satisfies relations \eqref{eq.KP.5.2} and \eqref{eq.KP.7.1}. It is easy to see that all the components of  ${\bf X}$ follow distributions $F_i \in \mathcal{D} \cap \mathcal{P_D}$, hence $0<\beta_{F_i} \leq \alpha_{F_i} < \infty$, for any $i = 1,\,\ldots,\,n$. We notice that while class $\mathcal{D}_n$ has the property of $\mathcal{D}$ for the linear combinations of its components, see Remark \ref{rem.KP.5.1},  however it is not known if this is true for the classes $\mathcal{P_D}_n$ and $\mathcal{D} \cap \mathcal{P_D}$, even in the case when the linear combinations are degenerated to unity. The reason is that very few papers up to now considered the closure properties in class  $\mathcal{P_D}_n$ with respect to sum, and these are restricted to the independence case, see \cite{bardoutsos:konstantinides:2011}. In the next result we provide a closure property with respect to the scalar product for the class $\mathcal{D}_n \cap \mathcal{P_D}_n $.

\bco \label{cor.KP.7.3}
Let ${\bf X}$  be random vector and $Y$ be random variable, with distributions ${\bf F}$ and $G$ respectively, and assume $G(0-)=0\,,\; G(0)<1$. Assume that ${\bf X}$ has identically distributed components, with common distribution $F$. If ${\bf X}$ and $Y$ are independent, with ${\bf F} \in \mathcal{D}_n \cap \mathcal{P_D}_n$ satisfying Assumption \ref{ass.KP.5.1}, then it holds $
{\bf H} \in \mathcal{D}_n \cap \mathcal{P_D}_n$. 
\eco

\pr~
It follows from Theorem \ref{th.KP.5.1} and Theorem \ref{th.KP.7.1}.
~\halmos

\bre \label{rem.KP.6.2}
 In several applications of risk theory, can be found the following assumption $\nu(({\bf 1},\,(\infty,\,\ldots,\,\infty)])> 0$, for the Radon measure $\nu$ of $MRV$, see for example \cite{konstantinides:li}, \cite{cheng:konstantinides:wang:2024} among others. As was mentioned in \cite{konstantinides:li} or \cite[Sec. 2]{li:2022(b)}, this assumption follows directly from relation
\beam \label{eq.KP.a}
\lim_{\xto} \dfrac 1{\bF(x)} \PP[X_1>t_1\,x,\,\ldots,\,X_n>t_n\,x] = \nu(({\bf t},\,(\infty,\,\ldots,\,\infty)])>0\,, 
\eeam
for any ${\bf t} \in [0,\,\infty]^n\setminus \{{\bf 0}\}$, recalling that ${\bf X}\in MRV(\alpha,\,F,\,\nu)$. It is well known that relation \eqref{eq.KP.a} means asymptotically dependent components, which further implies that the marginals of $MRV$ have weak equivalent tails, namely $\bF_i(x) \asymp \bF_j(x)$, as $\xto$, for any $i,\,j = 1,\,\ldots,\,n$, with $i\neq j$.

Let us denote 
\beao
MRV(F,\,\nu):= \bigcup_{0<\alpha< \infty} MRV(\alpha,\,F,\,\nu)\,,
\eeao 
and in case of asymptotic dependence we find that $MRV(F,\,\nu)$ belongs to $\mathcal{D}_n$, as also under a bit harder restriction, it belongs to $\mathcal{D}_n \cap \mathcal{P_D}_n$.
\ere

\bth \label{th.KP.6.2}
Let ${\bf X}=(X_1,\,\ldots,\,X_n)$ be a random vector with distribution ${\bf F} \in MRV(F,\,\nu)$ and its Radon measure $\nu$ be such that $\nu(({\bf 1},\,(\infty,\,\ldots,\,\infty)])>0$. Then it is true that
\begin{enumerate}
\item
${\bf F} \in \mathcal{D}_n$, hence 
\beam \label{eq.KP.6.b}
MRV(F,\,\nu) \subsetneq \mathcal{D}_n\,.
\eeam
\item
If distribution ${\bf F}$ has marginals, with strong equivalent tails, namely $\bF_i(x) \sim c_{ij}\,\bF_j(x)$, as $\xto$, with constants $c_{ij}>0$, for any $i,\,j = 1,\,\ldots,\,n$, then ${\bf F} \in \mathcal{P_D}_n$, that implies 
\beam \label{eq.KP.6.bc}
MRV(F,\,\nu) \subsetneq \mathcal{P_D}_n\,.
\eeam
\item
If the assumptions of (2) hold, then ${\bf F} \in \mathcal{D}_n \cap\mathcal{P_D}_n$, that implies $MRV(F,\,\nu) \subsetneq \mathcal{D}_n \cap\mathcal{P_D}_n$.
\end{enumerate}
\ethe

\pr~
\begin{enumerate}
\item
Let ${\bf F} \in MRV(\alpha,\,F,\,\nu)$ for some $\alpha\in (0,\,\infty)$. From relation \eqref{eq.KP.a}, we obtain  the inequality
\beam \label{eq.KP.b}
\limsup_{\xto} \dfrac{\PP[{\bf X} > {\bf b}\,{\bf t}\,x]}{\PP[{\bf X} > {\bf t}\,x]}=\limsup_{\xto} \dfrac{\nu[{\bf b}\,{\bf t},\,(\infty,\,\ldots,\,\infty)]\,\bF(x)}{\nu[{\bf t},\,(\infty,\,\ldots,\,\infty)]\,\bF(x)}\leq \widehat{b}^{-\alpha}< \infty\,,
\eeam
for any ${\bf b} \in (0,\,1)^n$, where $\widehat{b} :=\wedge_{i=1}^n b_i$, and in last step we used the homogeneity of the Radon measure $\nu$. Therefore, because of \eqref{eq.KP.b} and the fact that all marginals have weak equivalent tails we obtain $MRV (\alpha,\, F,\,\nu) \subsetneq \mathcal{D}_n$. Indeed, the marginals tails are weak equivalent because of the asymptotic dependence, which follows by relation $\nu(({\bf 1},\,(\infty,\,\ldots,\,\infty)])> 0$, see Remark \ref{rem.KP.6.2}, and the relation 
\beao
\mathcal{R}:=\bigcup_{0<\alpha <\infty} \mathcal{R}_{-\alpha} \subsetneq \mathcal{D}\,,
\eeao 
that implies any marginal can play the role of weak kernel. From the arbitrariness of $\alpha\in (0,\,\infty)$, we find \eqref{eq.KP.6.b}.

\item
Let ${\bf F} \in MRV(\alpha,\,F,\,\nu)$ for some $\alpha\in (0,\,\infty)$,  then, by relation \eqref{eq.KP.a} we obtain 
\beam \label{eq.KP.c} 
\limsup_{\xto}\dfrac{\PP[{\bf X}>{\bf v}\,{\bf t}\,x]}{\PP[{\bf X}>{\bf t}\,x]}=\limsup_{\xto} \dfrac{\nu[{\bf v}\,{\bf t},\,(\infty,\,\ldots,\,\infty)]\,\bF(x)}{\nu[{\bf t},\,(\infty,\,\ldots,\,\infty)]\,\bF(x)} \leq \widehat{v}^{-\alpha} < 1\,,
\eeam
where $\widehat{v} :=\wedge_{i=1}^n v_i >1$, for any ${\bf v} > {\bf 1}$. From the assumption that all marginals of ${\bf F}$ have strong equivalent tails and from $\mathcal{R} \subsetneq \mathcal{P_D}$, we find that any marginal can play the role of strong kernel of ${\bf F}$, which in combination with relation \eqref{eq.KP.c} gives that ${\bf F} \in \mathcal{P_D}_n$. Indeed, relation $\mathcal{R} \subsetneq \mathcal{P_D}$ can be verified from the uniqueness of the characterization of the classes by their Matuszewska indexes. Hence, $MRV(\alpha,\,F,\,\nu) \subsetneq \mathcal{P_D}_n$, which in combination with the arbitrarily chosen $\alpha\in (0,\,\infty)$ we conclude the \eqref{eq.KP.6.bc}.

\item
This follows immediately by combination of (1) with (2).~\halmos
\end{enumerate}

%% The Appendices part is started with the command \appendix;
%% appendix sections are then done as normal sections
%% \appendix

%% \section{}
%% \label{}

%% If you have bibdatabase file and want bibtex to generate the
%% bibitems, please use
%%
%%  \bibliographystyle{elsarticle-harv} 
%%  \bibliography{<your bibdatabase>}

%% else use the following coding to input the bibitems directly in the
%% TeX file.

\noindent \textbf{Acknowledgments.} 
The authors are grateful to prof. R. Leipus for several comments, that significantly improved the paper.

\end{document}